\documentclass[preprint,12pt]{elsarticle}
\usepackage[a4paper, total={6.1in, 9.5in}]{geometry}
\usepackage{amsmath,amssymb,amsthm}
\usepackage{color}
\usepackage{graphicx}
\usepackage{verbatim}
\usepackage{subfig}
\usepackage{multirow}
\usepackage{color}
\usepackage[rgb,dvipsnames]{xcolor}
\usepackage{mathtools}

\newcommand{\defeq}{\vcentcolon=}

\newcommand{\figdir}{mlffig}

\newtheorem*{remark}{Remark}

\newcommand{\cD}{{}^c\!D}
\newcommand{\cd}{{}^c\!\partial}
\newcommand{\strip}{\Omega}

\makeatletter
\def\ps@pprintTitle{%
 \let\@oddhead\@empty
 \let\@evenhead\@empty
 \def\@oddfoot{}%
 \let\@evenfoot\@oddfoot}
\makeatother

\begin{document}

\begin{frontmatter}
\title {Rational  Approximations for Oscillatory \\ Two-Parameter Mittag-Leffler Function}
\author[kfupm,Honain]{Aljowhara H. Honain}
\author[kfupm,kmfurati]{Khaled M. Furati}
\author[kfupm,ibrahim.sarumi]{Ibrahim O. Sarumi}
\author[mtsu]{Abdul Q. M. Khaliq}
\address[kfupm]{King Fahd University of Petroleum \& Minerals \\ 
Department of Mathematics \\ Dhahran, Saudi Arabia}
\address[Honain]{g201907570@kfupm.edu.sa}
\address[kmfurati]{kmfurati@kfupm.edu.sa}
\address[ibrahim.sarumi]{ibrahim.sarumi@kfupm.edu.sa}
\address[mtsu]{Middle Tennessee State University\\ Department of Mathematical Sciences\\ Murfreesboro, TN, USA\\ abdul.khaliq@mtsu.edu}
				
\begin{abstract}
The two-parameter Mittag-Leffler function $E_{\alpha, \beta}$ is of fundamental importance in fractional calculus. 
It appears frequently in the solutions of fractional differential and integral equations.
Nonetheless, this vital function is often expensive to compute. 
Several attempts have been made to construct cost-effective and accurate approximations. 
These attempts focus mainly on the completely monotone Mittag-Leffler functions. 
However, when $\alpha > 1$ the monotonicity property is largely lost and as such  roots and oscillations are exhibited.
Consequently, existing approximants constructed mainly for $\alpha \in (0,1)$ often fail to capture this oscillatory behavior.
In this paper, we construct computationally efficient and accurate rational approximants for $E_{\alpha, \beta}(-t)$, $t \ge 0$, with $\alpha \in (1,2)$.
This construction is fundamentally based on the decomposition of Mittag-Leffler function with real roots into one without and a polynomial.
Following which new approximants are constructed by combining the global Pad\'e approximation with a polynomial of appropriate degree.  
The rational approximants are extended to approximation of matrix Mittag-Leffler and different approaches to achieve efficient implementation for matrix arguments are discussed.
Numerical experiments are provided to illustrate the significant accuracy improvement achieved by the proposed approximants. 
\end{abstract}		
\end{frontmatter}	

\section{Introduction}

We are concerned with the approximation of the two parameter Mittag-Leffler function (MLF) $E_{\alpha,\beta}$, defined by
\begin{equation}\label{eq:111}
E_{\alpha,\beta}(z) =
\sum_{k=0}^{\infty} \frac{z^k}{\Gamma(\alpha k+\beta)}, \quad
\operatorname{Re} \alpha>0, \beta \in \mathbb{C}, \quad z \in \mathbb{C}.
\end{equation}
This entire function generalizes the MLF of one-parameter, $E_\alpha=E_{\alpha, 1}$, and contains several well-known 
special functions as special cases. In particular, $E_{1}$ is the exponential function, $E_{2}(-z^2)$ and $zE_{2,2}(-z^2)$ 
are the cosine and sine functions, respectively, among others. For some surveys on the Mittag-Leffler functions, see for example \cite{Kilbas2006,Gorenflo2014}. 

The Mittag-Leffler function $E_{\alpha, \beta}(z)$ arises frequently in the solutions of many physical problems described by differential and/or integral equations of fractional order. 
In the case $\alpha \in (1,2)$, the MLF appears naturally in the 
solutions of fractional diffusion-wave equations, 
fractional differential equation for motion, 
and fractional plasma equations, see \cite{gorenflo1996fractional,mainardi1996fractional,achar2001dynamics,stanislavsky2004fractional,stanislavsky2005twist}. 

Computing the MLF $E_{\alpha, \beta}$ is usually challenging 
(devising and implementing suitable algorithms) and expensive (computation time). Although the series \eqref{eq:111} converges for all $z \in \mathbb{C}$, it is impractical or ineffective to use it computationally for $|z| \geq 1$ because the series converges very slowly. 
As a result, various methods have been developed to evaluate the MLF. 
In \cite{Gorenflo2002}, an algorithm based on the location of the argument $z$ in the complex plane is developed, whereby for large $|z|$ values the asymptotic series as $|z| \to \infty$ is used, 
for $|z| < 1$ the series definition \eqref{eq:111} is used, and for values in intermediate regions the integral representation is used. 
In \cite{Garrappa2015}, an approach based on numerical inversion 
of Laplace transform is proposed. 
However, these approaches are often computationally time consuming, 
see \cite{Sarumi2020,Sarumi2021} for some CPU-time comparisons. 

Various rational approximations have been sought to approximate $E_{\alpha,\beta}$ 
efficiently and accurately.
Some of them are based solely on the series definition and provide accurate approximations 
for small argument values, see for example \cite{Iyiola2018a,Starovoitov2007,Borhanifar2015}.
Some others are based on global Pad\'e approach \cite{Winitzki2003} in which a hybrid of the local series 
definition and the asymptotic series representation is used.
These global approximations lead to approximants that are accurate over a wide range of arguments, 
see \cite{Atkinson2011,Zeng2015,Ingo2017,Sarumi2020,Sarumi2021}.

The existing rational approximants are effective when the MLF is completely monotone, which is the case when $\alpha \in (0,1)$ and $\beta \ge \alpha$. 
However, when $\alpha > 1$, as discussed in \cite{hanneken2007enumeration,hanneken2013alpha,duan2014zeros},
for some $\alpha$-$\beta$ combinations, $E_{\alpha,\beta}$ 
is often oscillatory and having multiple real zeros.
In these situations, the rational approximants might fail to trace the oscillation profile and match the large zeros.

In this paper, we consider in the $\alpha\beta$-plane the strip 
\begin{equation}
\label{eq:strip} 
\strip = \{(\alpha,\beta): 1<\alpha < 2, \beta \ge 1\}.
\end{equation} 
The regions where $E_{\alpha,\beta}$ is monotone and where it is oscillatory are classified by studying the zeros of its derivative.
In addition, we introduce new rational approximants that sufficiently capture more of the MLF roots and are able to trace their oscillations.
This is achieved by decomposing $E_{\alpha,\beta}$ into one without and a polynomial. 
Based on this decomposition, we construct new rational approximants
as a sum of the global Pad\'e approximant and a polynomial. 
Using this approach, we are able to capture sufficiently many roots of $E_{\alpha,\beta}$ by choosing the appropriate degree of the polynomial.
Further, we generalize these approximants to matrix MLF and discuss different approaches in which the approximant can be implemented for matrix arguments. 
Apart from the computation cost, all approaches yield close approximations.

Numerical experiments are presented to illustrate the performance of our approximants.
Considering MLF with scalar arguments, we observe the limitation of the existing global Pad\'e approximants in approximating the oscillatory MLFs. 
It is then demonstrated that the new approximants introduced here are able to improve 
the situation thereby providing sufficiently accurate approximations. 
For the matrix MLF, we consider four approaches for computing the rational approximants and compare them to the matrix MLF algorithm developed 
in \cite{Garrappa2018}.
Using the values obtained from this algorithm as a reference, it is seen that our 
approximants provide close values with significantly less computation time.

The rest of the paper is organized as follows: a characterization of the oscillatory behavior of the MLF 
is given in Section \ref{sec:charaterization}. In Section \ref{sec:derooting}, 
the derooting representation is discussed. 
Global Pad\'e approximants of MLF with scalar and matrix arguments 
are discussed in sections \ref{sec:rational approx} and \ref{sec:matrix MLF}, respectively. 
Some applications to the evaluation of solutions of fractional 
evolution equations are discussed in Section \ref{sec:applications}.

\section{Monotonicity and oscillatory properties}
\label{sec:charaterization}

In this section, we characterize the oscillatory behavior of 
$$
E_{\alpha,\beta }(-t), \qquad (\alpha,\beta) \in \strip, \qquad t \ge 0.
$$
For the discussion that follows, the following formula for the derivative \cite{Kilbas2006} plays a fundamental role:
\begin{equation} 
\label{eq:deriv}
\frac{d}{d t}E_{\alpha, \beta}(-t)=-\frac{1}{\alpha}\left[E_{\alpha, \beta+\alpha-1}(-t)-(\beta-1) E_{\alpha, \beta+\alpha}(-t)\right].
\end{equation}
Since $E_{\alpha,\beta}(0)= 1 / \Gamma(\beta)$, it follow from \eqref{eq:deriv} that
\begin{equation}
\label{eq:deriv0}
\frac{d}{d t}E_{\alpha, \beta}(0)=
-\frac{1}{\alpha}
\left[\frac{1}{\Gamma(\beta+\alpha-1)} -
\frac{\beta-1}{ \Gamma(\beta+\alpha)}  \right] = -\frac{1}{\Gamma(\beta+\alpha)} < 0.
\end{equation}

The function $E_{\alpha,\beta}(-t)$, $t \ge 0$, is known to be completely monotone for $0 < \alpha \le 1$, $\beta \ge \alpha$. 
However, this monotonicity is not necessarily preserved when $(\alpha,\beta)\in \strip$ since either the function or its derivative could have roots.
On the other hand, it has been established in \cite{hanneken2013alpha} that $E_{\alpha,\beta}$ has at most a finite number of roots when $(\alpha,\beta)\in \strip$.
Thus, $E_{\alpha,\beta}(-t)$ is monotone for sufficiently large $t$.
Furthermore, it is shown in \cite{hanneken2013alpha} that 
there exists only one zero when $\alpha$ is sufficiently close to $1$ 
while the number of zeros increases as $\alpha$ increases towards 2.

\begin{figure}[p]
\centering
\includegraphics[width=.7\textwidth]{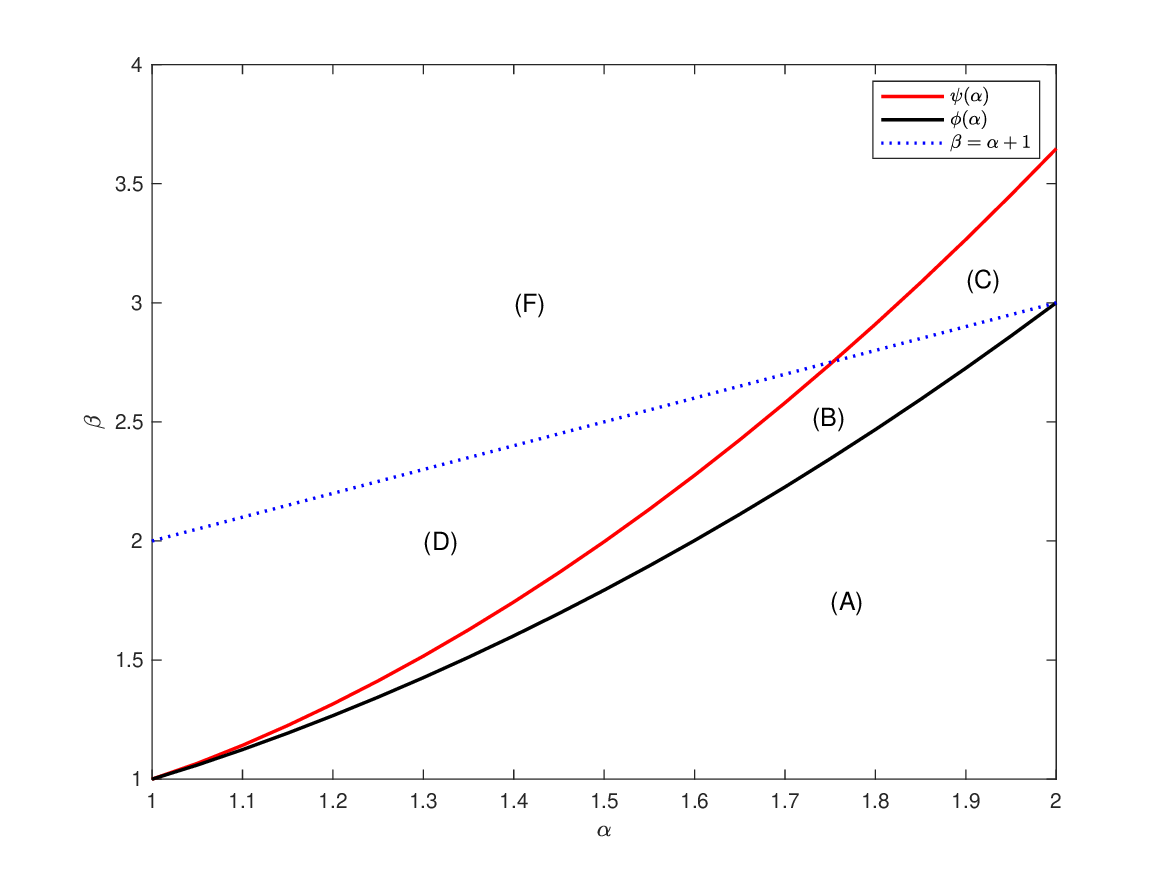}
\caption{$\alpha$-$\beta$ phase diagram for $E_{\alpha, \beta}(-t)$ and its derivative.}
\label{fig:diag}
\end{figure}

In Figure \ref{fig:diag}, the curve $\phi$ is the boundary given in
\cite[Table 1]{hanneken2013alpha} such that $E_{\alpha,\beta}$ has a finite number of real roots when $(\alpha,\beta)$ is below it and none above it.
However, the function in the region below $\phi (\alpha)$ differs with respect to the number of roots. 
Typically, the number of roots increases as $\alpha$ goes to 2. 
Similarly, we constructed the boundary $\psi$ for the derivative.
The corresponding data is presented in Table \ref{tab:beta}.

The monotonicity and oscillatory behavior of $E_{\alpha,\beta}$ is characterized by the different regions 
marked in Figure \ref{fig:diag} and described in Table \ref{table:nonlin}.
It follows from \eqref{eq:deriv0} that $E_{\alpha,\beta}(-t)$ is monotonically decreasing for $(\alpha,\beta)$ in the regions (D) and (F). 
In the regions (B) and (C), although the function has no real roots, it could have a finite number of oscillations due to the roots of the derivative.

\begin{remark}
Since $E_{\alpha,\beta}(-t) \to 0$ as $t \to \infty$, it is obvious that each root of
it is followed by at least one root of its derivative.
This is consistent with the inequality
$$
\phi (\alpha) < \psi (\alpha). 
$$ 
\end{remark}

\begin{table}[p]
\begin{center}
\caption{The data for the boundaries $\phi$ and $\psi$.}
\label{tab:beta}
\begin{tabular}{l@{\hskip 1in}c@{\hskip 1in}c}
	\hline
	$\alpha$ & $\phi(\alpha)$ & $\psi(\alpha)$  \\
	\hline
	$1.00$   & $ 1.00000  $&$1.00000$\\
	$1.05$ & $ 1.05924  $&$1.06640$  \\
	$1.10$ & $ 1.12400  $& $1.14204$ \\
	$1.15$ & $  1.19325 $& $  1.22532 $  \\
	$1.20$ & $ 1.26674  $& $  1.31565 $  \\
	$1.25$ & $  1.34437 $& $ 1.41277  $ \\
	$1.30$ & $1.42608   $& $ 1.51654  $  \\
	$1.35$ & $1.51187   $& $ 1.62689  $  \\
	$1.40$ & $  1.60173 $& $  1.74374 $  \\
	$1.45$ & $ 1.69565  $& $ 1.86706  $  \\
	$1.50$ & $  1.79365 $& $   1.99685$  \\
	$1.55$ & $  1.89573 $& $  2.13306 $  \\
	$1.60$ & $   2.00191$& $ 2.27568  $  \\
	$1.65$ & $  2.11219 $& $  2.42471  $  \\
	$1.70$ & $  2.22660 $& $2.58014   $ \\
	$1.75$ & $   2.34513$& $ 2.74196 $  \\
	$1.80$ & $ 2.46779  $& $ 2.91017  $ \\
	$1.85$ & $  2.59460 $& $  3.08477 $  \\
	$1.90$ & $  2.72557 $& $  3.26575 $  \\
	$1.95$ & $  2.86070 $& $  3.45311 $  \\
	$2.00$& $  3.00000 $ & $  3.64686 $  \\
	\hline
\end{tabular}
\end{center}
\end{table}

\begin{table}
\caption{Number of roots in the regions of Figure \ref{fig:diag}.}
\centering
\begin{tabular}{c c c}
Region & $E_{\alpha, \beta}(-t)$ & $d E_{\alpha, \beta}/dt$ \\
\hline
A & Finite number of roots & Finite number of roots  \\
B & No real roots & Finite number of roots  \\
C & No real roots & Finite number of roots   \\
D & No real roots & No real roots \\ 
F & No real roots  & No real roots     \\ [1ex]
\hline
\end{tabular}
\label{table:nonlin}
\end{table}

For the region below the line $\beta = \alpha +1$, which includes the regions (A), (B) and (D), the following decomposition is presented in  \cite{hanneken2007enumeration},
\begin{equation}	
\label{eq:decfg}
E_{\alpha, \beta}(-t)=g_{\alpha, \beta}(-t)+f_{\alpha, \beta}(-t),
\qquad t > 0, \quad \beta < \alpha + 1,
\end{equation}
where the function $f_{\alpha, \beta}(-t)$ is asymptotically approaching zero as 
$t \rightarrow \infty$ and $g_{\alpha, \beta}(-t)$ is an oscillatory function.

When $\beta = 1$ and $1<\alpha<2$, the function $f_{\alpha,1}(-t)$  
is known to be a completely monotonic function while the function $g_{\alpha,1}(-t)$ is oscillatory with exponentially decreasing amplitude \cite{gorenflo1996fractional}.
The properties of $E_{\alpha}$, $1< \alpha < 2$, such as monotonicity, roots, and the oscillatory behavior have been discussed in details in \cite{honain2023rational}. 

\begin{figure}
\centering
\includegraphics[width=0.5\textwidth]{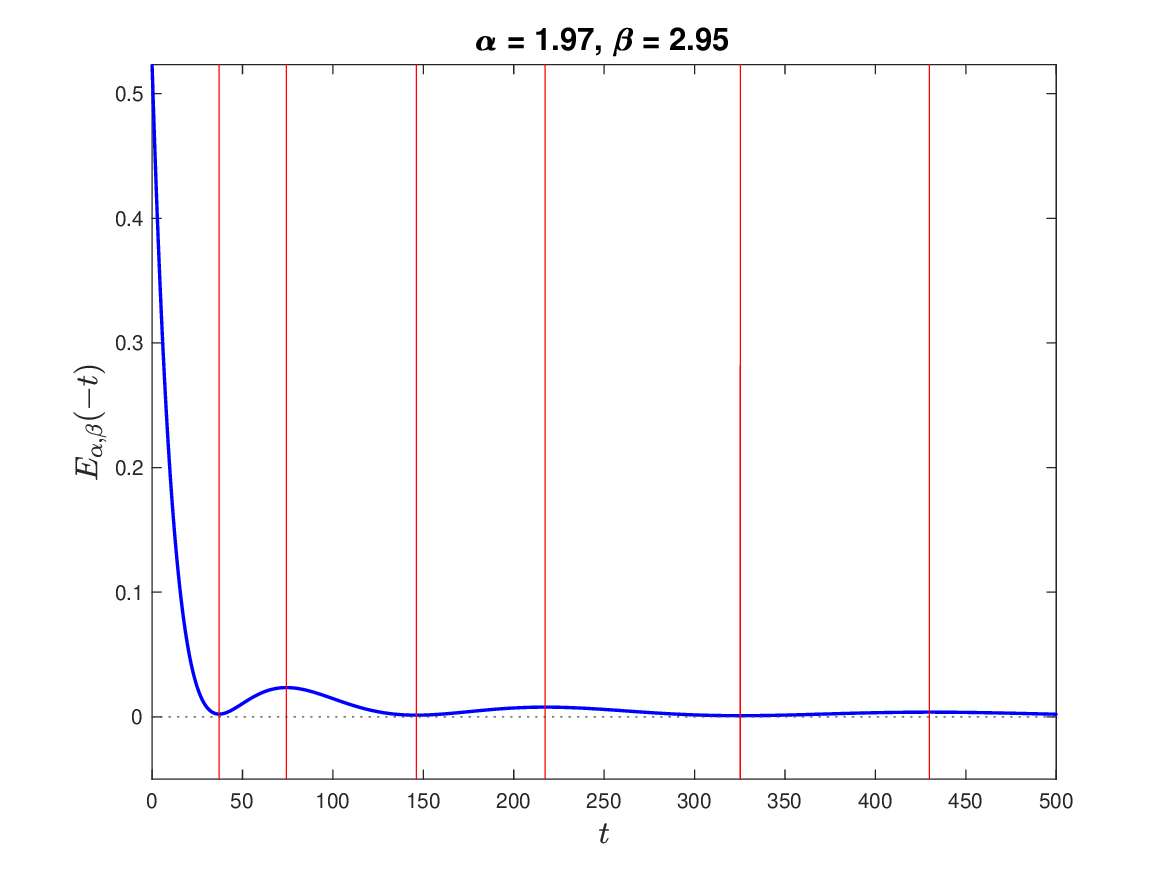}
\caption{Typical oscillations in region (B) about the roots (marked by the red lines) of the derivative.}
\label{fig:Bregion}
\end{figure}

In the region (B), although $E_{\alpha,\beta}(-t)$ has no roots as it decays to zero, it undergoes oscillations due to the roots of its derivative as demonstrated in
Figure \ref{fig:Bregion}.
This behavior is also consistent with the oscillatory nature of $g_{\alpha,\beta}$ there.
On the other hand, in the region (D) in which $E_{\alpha,\beta}$ is monotone, the behavior of the oscillatory function $g_{\alpha, \beta}$ fades away and $f_{\alpha, \beta}(-t)$ dominates.

\section{Derooting decomposition} 
\label{sec:derooting}

A key step to accurately approximate $E_{\alpha,\beta}(-t)$ when 
$1 < \alpha < 2$ is to capture its oscillatory behavior and roots.
However, in region (A) the number of roots in general could exceed the 
ones of rational approximants. This makes these approximants valid only 
for small arguments.

A Derooting approach can be used to extend the validity of rational approximants to large intervals.
For this purpose, we use the following recursive identity \cite{haubold2011mittag}
\begin{equation}
	\label{eq:rec1}
	\begin{split}
		E_{\alpha, \beta}(t) &= t^r E_{\alpha,\beta+\alpha r}(t) +
		P_{\alpha, \beta}^{r-1}(t), 
		\\ & 
		\alpha>0, \quad \beta > 0, \quad t\in \mathbb{R}, \quad  r = 0, 1, 2, \dots,
	\end{split}
\end{equation}
where,
\begin{equation}\label{eq:mlf poly term}
	P_{\alpha,\beta}^{r-1}(t) :=
	\left\{
	\begin{array}{ll}
		\sum_{k=0}^{r-1} \dfrac{t^{k}}{\Gamma(\alpha k+\beta)},  
		& r=1,2,3, \dots\\
		0, & r=0.
	\end{array}
	\right.
\end{equation}  
In this identity, $E_{\alpha,\beta}$ is decomposed into another MLF
with shifted second parameter and a polynomial of degree $r-1$.
This allows us to replace the parameters $(\alpha,\beta)$ in region (A) by 
parameters $(\alpha,\beta+\alpha r)$ above the boundaries $\phi(\alpha)$ or 
above $\psi(\alpha)$.
Consequently, we obtain a function $E_{\alpha,\beta +\alpha r}$ that has no roots or even could be monotone.
Therefore,
the approximation of $E_{\alpha,\beta}(-t)$, $1<\alpha<2$, with a finite number of real roots can be replaced by approximating another one that has no roots or one that is monotone.

\section{Rational approximation}
\label{sec:rational approx}

Based on the global rational approximation technique introduced for transcendental 
functions in \cite{Winitzki2003}, a variety of global Pad\'e approximants
for $E_{\alpha,\beta} (-t)$, $0 < \alpha < 1$, are developed and implemented 
\cite{Atkinson2011,Zeng2015,Ingo2017, Sarumi2020, Sarumi2021}.
These approximants are still valid for $E_{\alpha,\beta} (-t)$, $1 < \alpha <2$, 
over intervals in which the number of its roots does not exceed 
the number of approximant roots.
This renders these approximants inadequate when approximating MLF over an
interval in which it has multiple roots, as the case when $(\alpha,\beta)$ in 
region (A)
of Figure \ref{fig:diag}.

In this section we consider the global Pad\'e approximants 
$R_{\alpha, \beta}^{7,2}$ and $R_{\alpha, \beta}^{13,4}$ developed in 
\cite{Sarumi2020, Sarumi2021} and examine their performance when $1 < \alpha < 2$.
Then we describe how the derooting decomposition leads to approximants with increased number of roots and thus able to track more of MLF roots.

\subsection{Global Pad\'e approximation}

For $\{(\alpha,\beta): 1<\alpha < 2, \beta \ge 1, \alpha \ne \beta\}$, 
we consider the global Pad\'e approximants 
\begin{equation}
\label{eq:R72}
R_{\alpha, \beta}^{7,2}(t)=\frac{1}{\Gamma(\beta-\alpha)}\frac{p_1+p_2 t+p_3 t^2 +t^3}{q_0+q_1 t+q_2 t^2+q_3 t^3+t^4}, \quad t \geq 0,
\end{equation}
and 
\begin{equation}
\label{eq:R134}
R_{\alpha, \beta}^{13,4}(t)=\frac{1}{\Gamma(\beta-\alpha)}\frac{p_1+p_2 t+\cdots+p_7 t^6+t^7}{q_0+q_1 t+q_2 t^2+\cdots+q_7 t^7+t^8}, \quad t \geq 0,
\end{equation}
where the coefficients satisfy respectively the systems
$$
\left(\begin{array}{lllllll}1 & 0 & 0 & a_0 & 0 & 0 & 0  
	\\ 
	0 & 1 & 0 & a_1 & a_0 &  0 & 0 
	\\ 
	0 & 0 & 1 & a_2 & a_1 & a_0 & 0  
	\\ 
	0 & 0 & 0 & a_3 & a_2 & a_1 & a_0 
	\\ 
	0 & 0 & 0 & a_4 & a_3 & a_2 & a_1 
	\\ 
	0 & 0 & 0 & a_5& a_4 & a_3 & a_2  
	\\ 
	0 & 0 & 1 & 0 & 0 & 0 &b_0 
\end{array}\right)\left(\begin{array}{l}p_1 \\ p_2 \\ p_3 \\ q_0 \\ q_1 \\ q_2 \\ q_3\end{array}\right)=\left(\begin{array}{c}0 \\ 0 \\ 0 \\  -1 \\ -a_0 \\ -a_1 \\ -b_1 \end{array}\right),
$$
$$
\left(\begin{array}{llllllllllllllllll}1 & 0 & 0 & 0 & 0 & 0 & 0 & a_0 & 0 & 0 & 0 & 0 & 0 & 0 & 0 \\ 0 & 1 & 0 & 0 & 0 & 0 & 0 & a_1 & a_0 & 0 & 0 & 0 & 0 & 0 & 0 \\ 0 & 0 & 1 & 0 & 0 & 0 & 0 & a_2 & a_1 & a_0 & 0 & 0 & 0 & 0 & 0 \\ 0 & 0 & 0 & 1 & 0 & 0 & 0 & a_3 & a_2 & a_1 & a_0 & 0 & 0 & 0 & 0 \\ 0 & 0 & 0 & 0 & 1 & 0 & 0 & a_4 & a_3 & a_2 & a_1 & a_0 & 0 & 0 & 0 \\ 0 & 0 & 0 & 0 & 0 & 1 & 0 & a_5 & a_4 & a_3 & a_2 & a_1 & a_0 & 0 & 0 \\ 0 & 0 & 0 & 0 & 0 & 0 & 1 & a_6 & a_5 & a_4 & a_3 & a_2 & a_1 & a_0 & 0 \\ 0 & 0 & 0 & 0 & 0 & 0 & 0 & a_7 & a_6 & a_5 & a_4 & a_3 & a_2 & a_1 & a_0 \\ 0 & 0 & 0 & 0 & 0 & 0 & 0 & a_8 & a_7 & a_6 & a_5 & a_4 & a_3 & a_2 & a_1 \\ 0 & 0 & 0 & 0 & 0 & 0 & 0 & a_9 & a_8 & a_7 & a_6 & a_5 & a_4 & a_3 & a_2 \\ 0 & 0 & 0 & 0 & 0 & 0 & 0 & a_{10} & a_9 & a_8 & a_7 & a_6 & a_5 & a_4 & a_3 \\ 0 & 0 & 0 & 0 & 0 & 0 & 0 & a_{11} & a_{10} & a_9 & a_8 & a_7 & a_6 & a_5 & a_4 \\ 0 & 0 & 0 & 0 & 1 & 0 & 0 & 0 & 0 & 0 & 0 & 0 & b_0 & b_1 & b_2 \\ 0 & 0 & 0 & 0 & 0 & 1 & 0 & 0 & 0 & 0 & 0 & 0 & 0 & b_0 & b_1 \\ 0 & 0 & 0 & 0 & 0 & 0 & 0 & 0 & 0 & 0 & 0 & 0 & 0 & 0 & b_0\end{array}\right)\left(\begin{array}{l}p_1 \\ p_2 \\ p_3 \\ p_4 \\ p_5 \\ p_6 \\ p_7 \\ q_0 \\ q_1 \\ q_2 \\ q_3 \\ q_4 \\ q_5 \\ q_6 \\ q_7\end{array}\right)=\left(\begin{array}{c}0 \\ 0 \\ 0 \\ 0 \\ 0 \\ 0 \\ 0 \\ -1 \\ -a_0 \\ -a_1 \\ -a_2 \\ -a_3 \\ -b_3 \\ -b_2 \\ -b_1\end{array}\right),
$$
$$
a_j=\frac{(-1)^{j+1} \Gamma(\beta-\alpha)}{\Gamma(\beta+j \alpha)}, \qquad b_j=\frac{(-1)^{j+1} \Gamma(\beta-\alpha)}{\Gamma(\beta-(j+1) \alpha)}, \,\, j=0,1, \ldots \, .
$$

Computationally, it is observed that $R_{\alpha,\beta}^{13,4}$ has at most three real roots, while $R_{\alpha,\beta}^{7,2}$ has exactly one real root. 
Consequently, in approximating some of the oscillatory MLFs, each approximant 
is only able to trace the function up to its largest root. 
To illustrate the accuracy and limitations of both approximants when $1<\alpha<2$, we compare with
the Matlab routine ml \cite{Garrappa2015} and ml\_matrix \cite{Garrappa2018}.

\begin{figure}[p]
\centering
\subfloat[Region(A): one root]
{\includegraphics[width=0.49\textwidth]{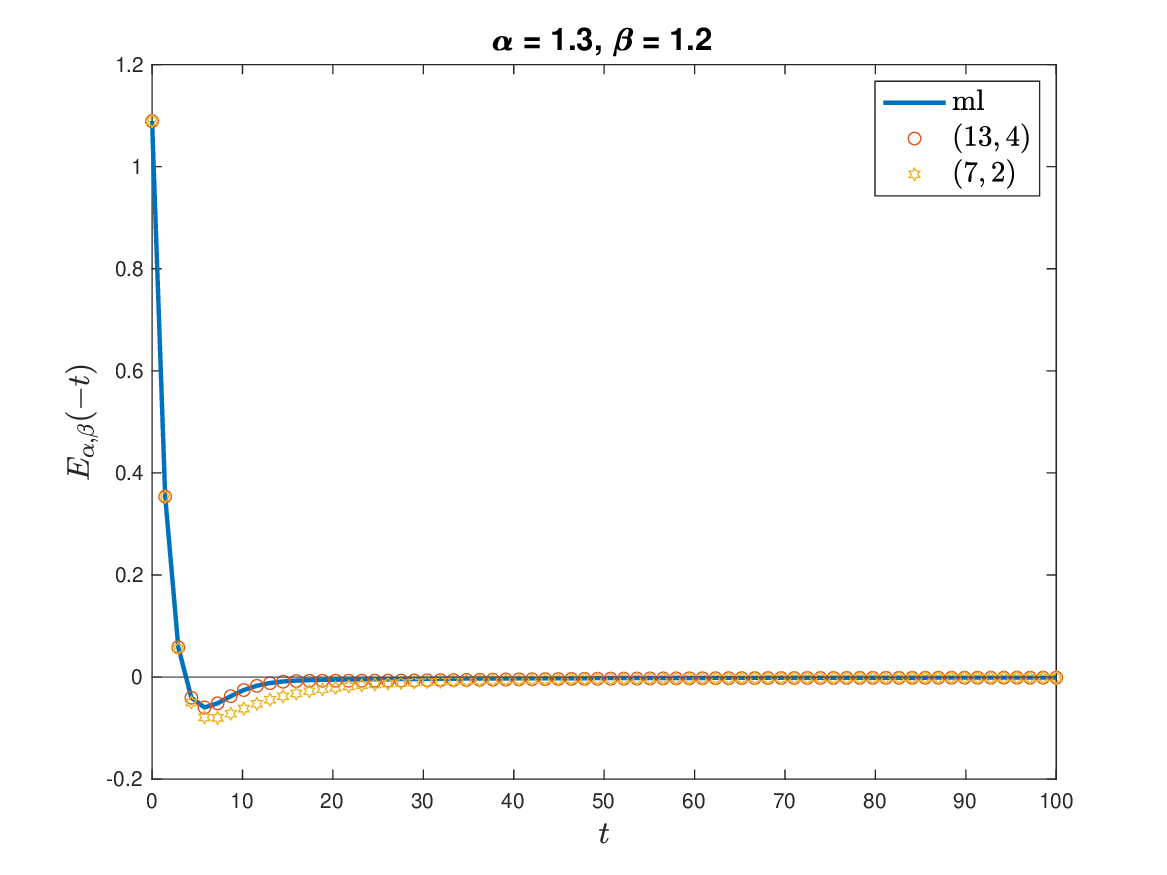}}
\subfloat[Region(A): many roots]
{\includegraphics[width=0.49\textwidth]{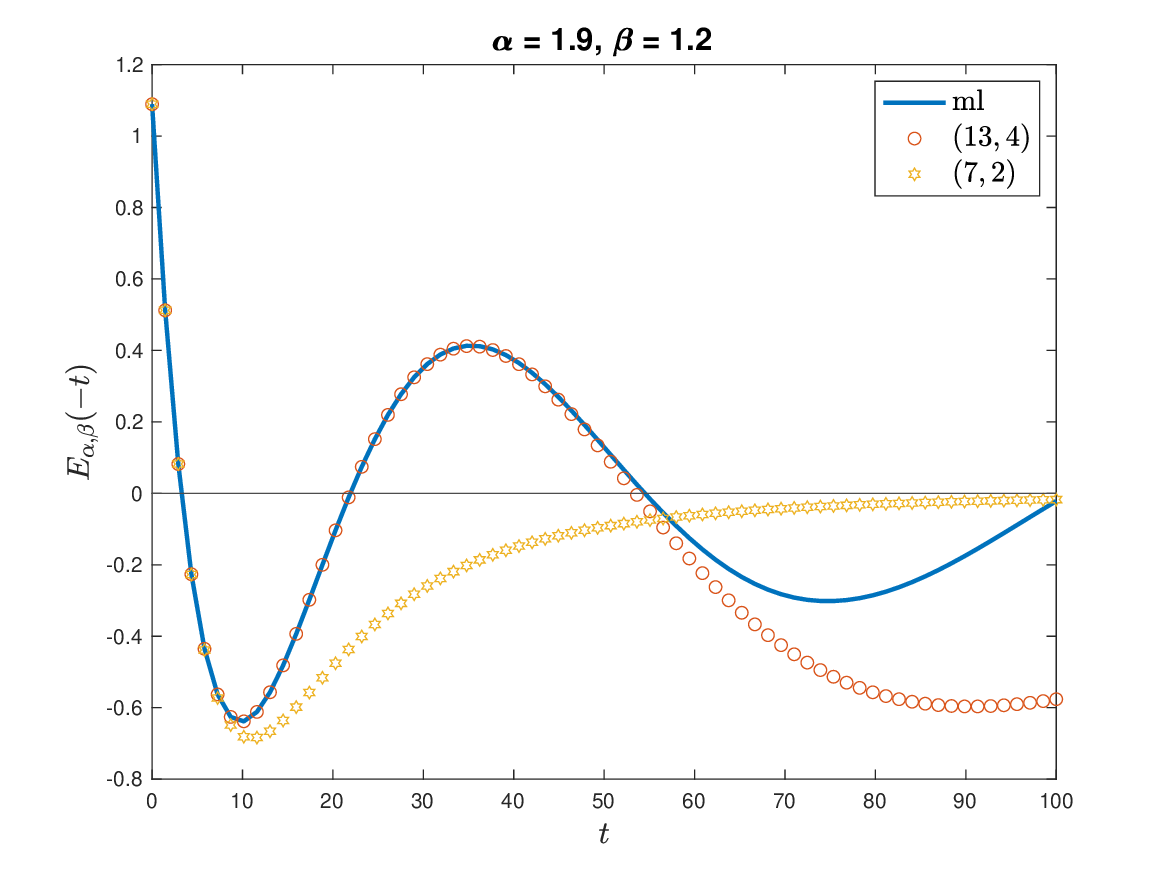}}
\caption{
Plots of $R_{\alpha,\beta}^{13,4}$ and $R_{\alpha,\beta}^{7,2}$ 
approximants of $E_{\alpha,\beta}$, $(\alpha,\beta)$ in region (A).}
\label{fig:1}
\end{figure}

In Figure \ref{fig:1}, two examples from region (A) are shown.
In plot (a), the MLF has only one root which is typically the case for $1<\alpha < 1.4$. As observed, $R_{\alpha,\beta}^{13,4}$ provides a good approximation in
this cases.
In plot (b), the MLF has many roots which typically the case as $\alpha$ 
gets closer to 2.
As expected, each approximant fails beyond its largest root to capture the oscillations of the MLF .

\begin{figure}[p]
\centering
\subfloat[Region(B): no roots and oscillatory]
{\includegraphics[width=0.49\textwidth]{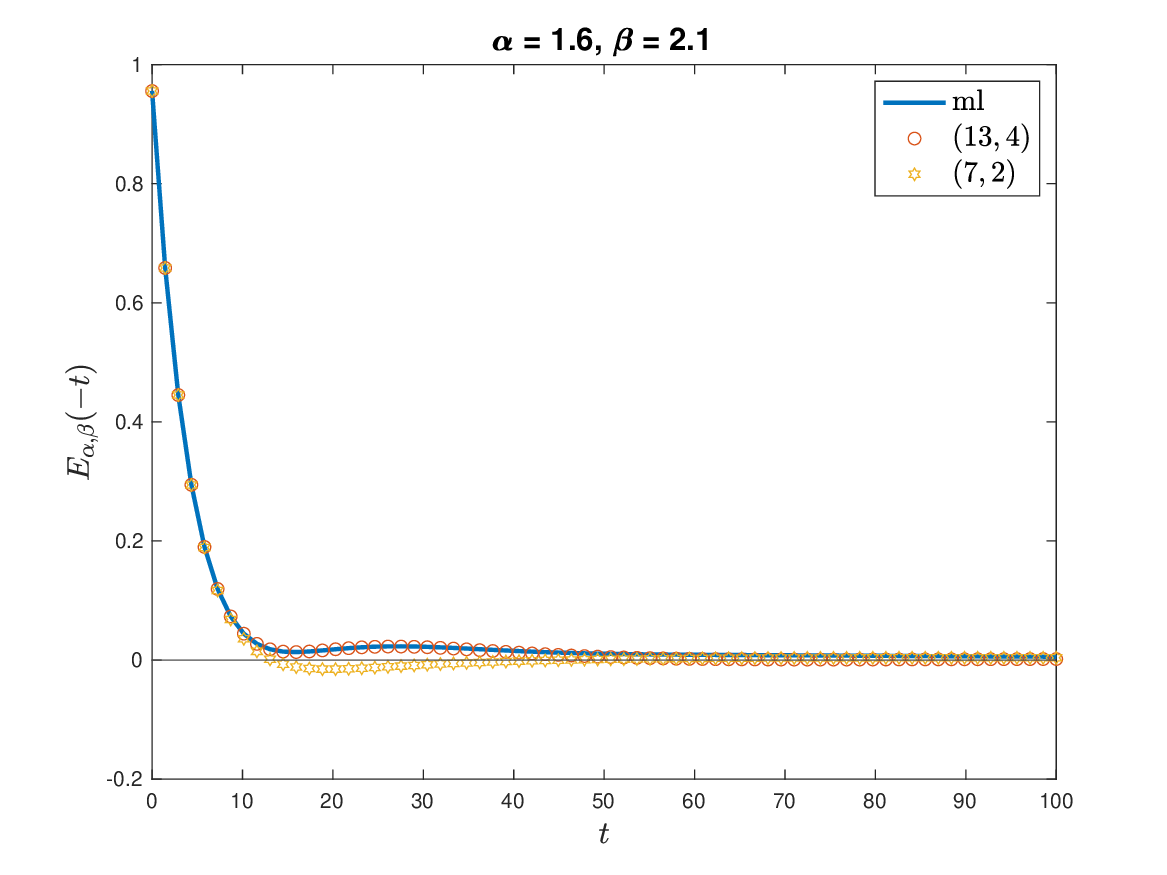} }
\subfloat[Region(C): no roots and oscillatory]
{\includegraphics[width=0.49\textwidth]{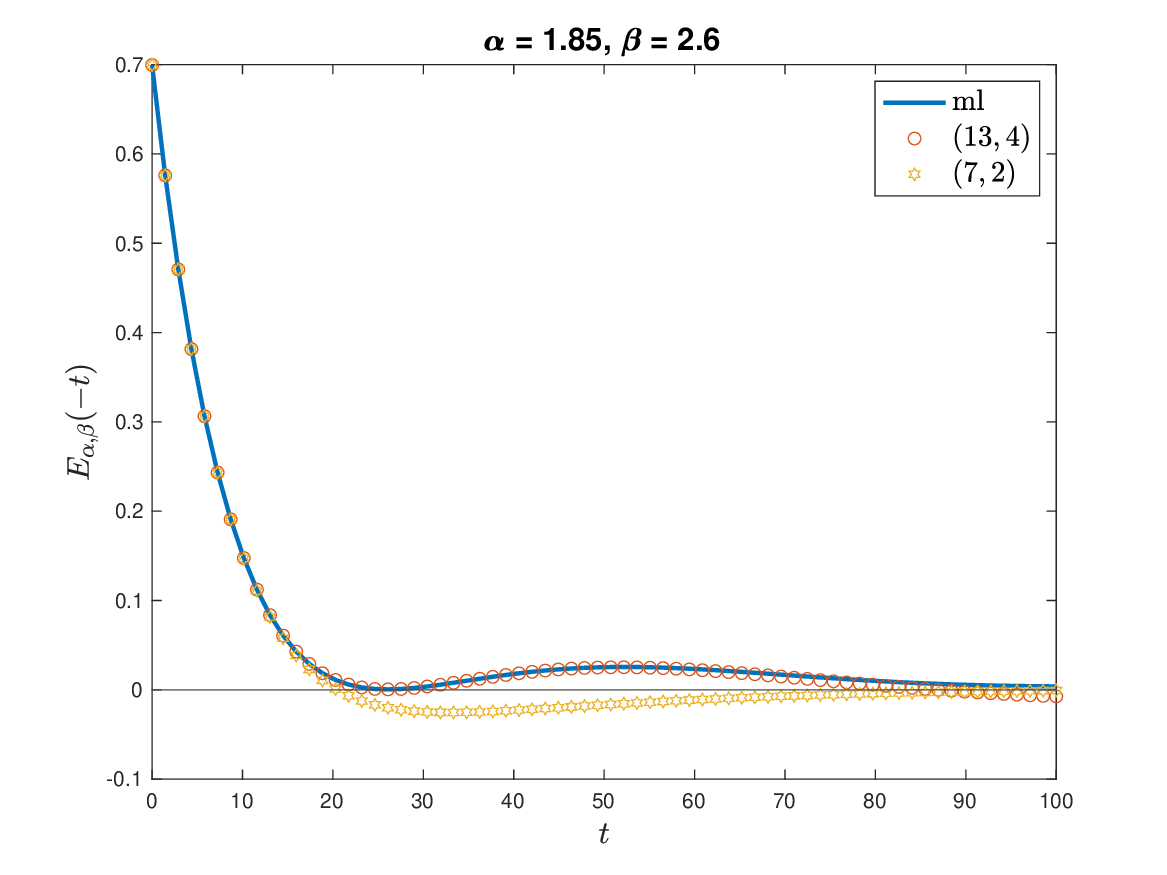} }
\caption{Plots of $R_{\alpha,\beta}^{13,4}$ and $R_{\alpha,\beta}^{7,2}$ 
approximants of $E_{\alpha,\beta}$, $(\alpha,\beta)$ in
regions (B) and (C).}
\label{fig:2}
\end{figure}

In Figure \ref{fig:2}, typical cases from Region (B) and (C) are shown.
The corresponding MLFs in both regions is root-free but are oscillatory. 
As shown, $R_{\alpha,\beta}^{13,4}$ is sufficiently accurate for an extended interval. However, eventually it takes negative values as observed in plot (b).

\begin{figure}[p]
\centering
\subfloat[Region(D): monotone and positive]
{\includegraphics[width=0.49\textwidth]{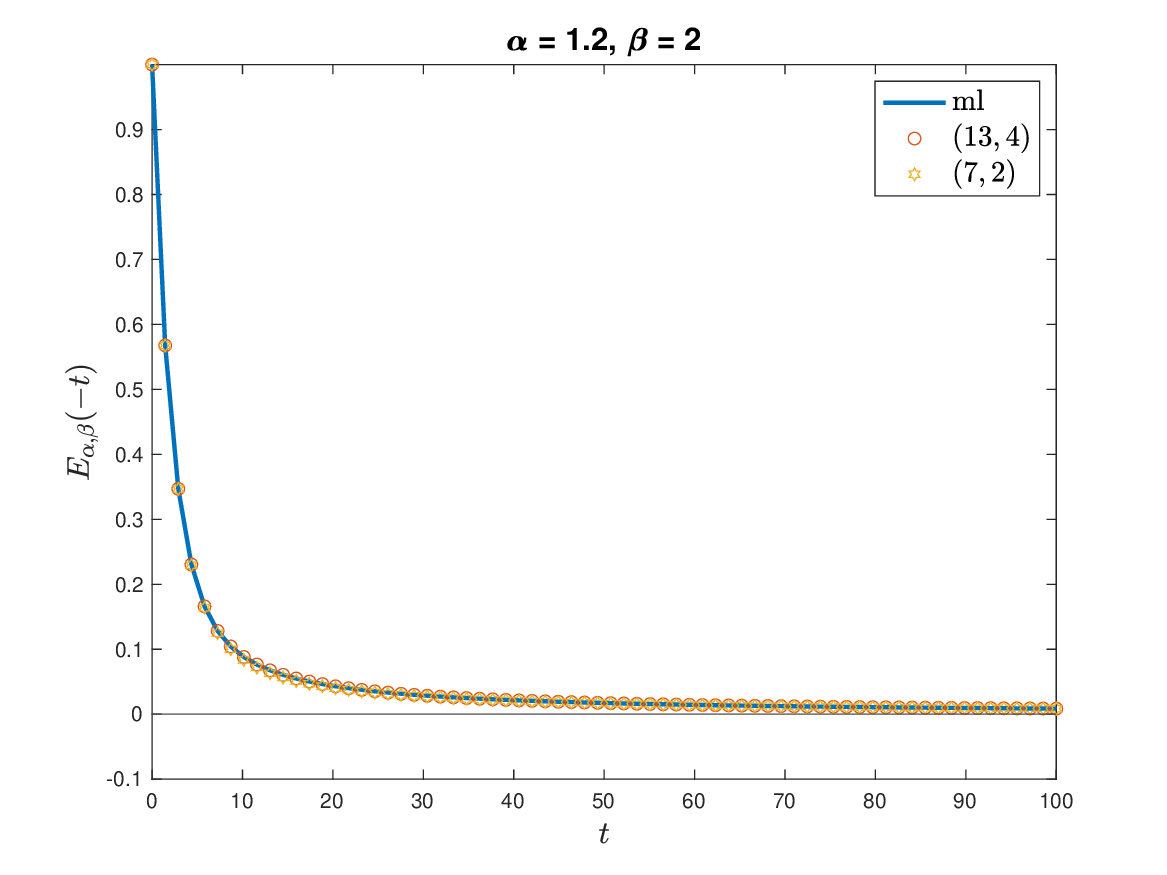} }
\subfloat[Region(F): monotone and positive]
{\includegraphics[width=0.49\textwidth]{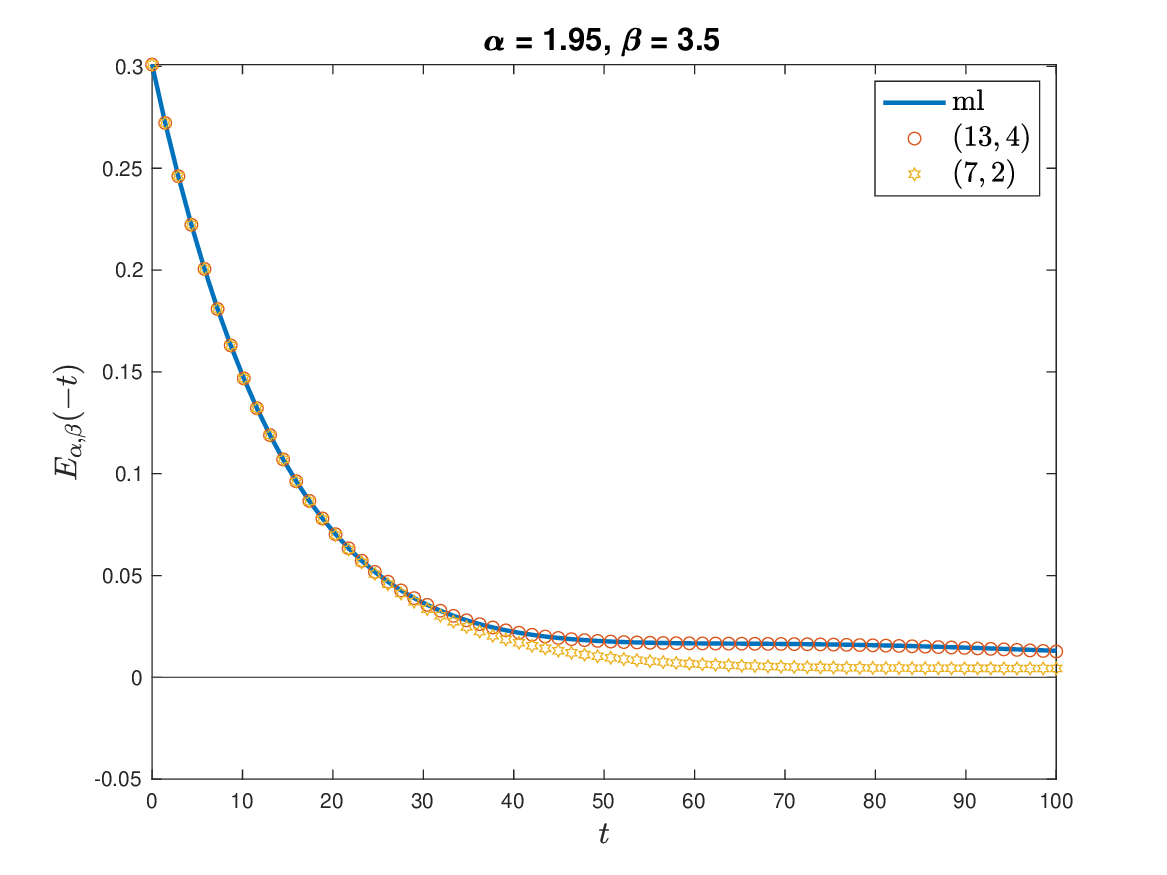} }
\caption{Plots of $R_{\alpha,\beta}^{13,4}$ and $R_{\alpha,\beta}^{7,2}$ 
approximants of $E_{\alpha,\beta}$, $(\alpha,\beta)$ in regions (D) and (F).}
\label{fig:3}
\end{figure}

In Figure \ref{fig:3}, typical cases from regions (D) and (F) are illustrated.  
In both regions the MLF is monotone and globally positive.
As seen clearly, $R_{\alpha,\beta}^{13,4}$ is a good approximation for
an extended interval.

As general guidelines, we have the following.
\begin{itemize}
\item 
When $(\alpha,\beta)$ in Regions (A), the existing global Pad\'e approximants 
are not effective over extended intervals due to the existence of roots.
\item 
When $(\alpha,\beta)$ in region (B) or (C), more accurate approximants should be 
used due to the oscillatory behavior of MLF.
\item 
When $(\alpha,\beta)$ in Region (D) or (F), the approximant $R_{\alpha,\beta}^{13,4}$
is sufficiently accurate since MLF is globally monotone.
\end{itemize}

\subsection{Rational approximation for oscillatory MLFs}

As observed in the above subsection, when $(\alpha,\beta)$ lies below the boundary $\psi(\alpha)$ 
in Figure \ref{fig:diag},
then more accurate and multi roots approximants are essential for extended intervals.
Next, we propose a class of rational approximants that have the 
desirable performance.

It follows from \eqref{eq:rec1} that we can write 
$$
E_{\alpha, \beta}(-t) = (-t)^r E_{\alpha,\beta+\alpha r}(-t) + P_{\alpha, \beta}^{r-1}(-t), \qquad t > 0.
$$
Let $R^{m,n}_{\alpha,\beta} $ be the global Pad\'e approximants described in \cite{Sarumi2020}.
We introduce the rational approximant 
\begin{equation}
\label{eq:approx}
 E_{\alpha, \beta}(-t) \approx R_{\alpha, \beta}^{m,n,r}(t), \qquad t > 0,
\end{equation}
where,  
\begin{equation}
\label{eq:general}
R_{\alpha,\beta}^{m,n,r}(t) := (-t)^r R_{\alpha, \beta + \alpha r}^{m,n}(t) + P_{\alpha, \beta}^{r-1}(-t), 
\quad r=0, 1,2, \dots.
\end{equation}
Note that as a special case we have
$$
R_{\alpha,\beta}^{m,n,0} = R^{m,n}_{\alpha,\beta}.
$$
The integer $r$ could be chosen so that $(\alpha,\beta+\alpha r)$ lies above the
boundary $\psi(\alpha)$ in Figure~\ref{fig:diag}.
As demonstrated in the previous subsection, this shifting of parameters would enhance the 
performance of $R^{m,n}_{\alpha,\beta+\alpha r}$ since in this case 
$E_{\alpha,\beta+\alpha r}(-t)$ is root-free and non-oscillatory.
Furthermore, $r$ could be chosen large enough so that the approximant 
$R_{\alpha, \beta}^{m,n,r}$ captures the desired number of roots and oscillations of $E_{\alpha,\beta}$ 
over an extended interval.

\begin{figure}[p]
\centering
\includegraphics[width=.49\textwidth,height=.3\textwidth]{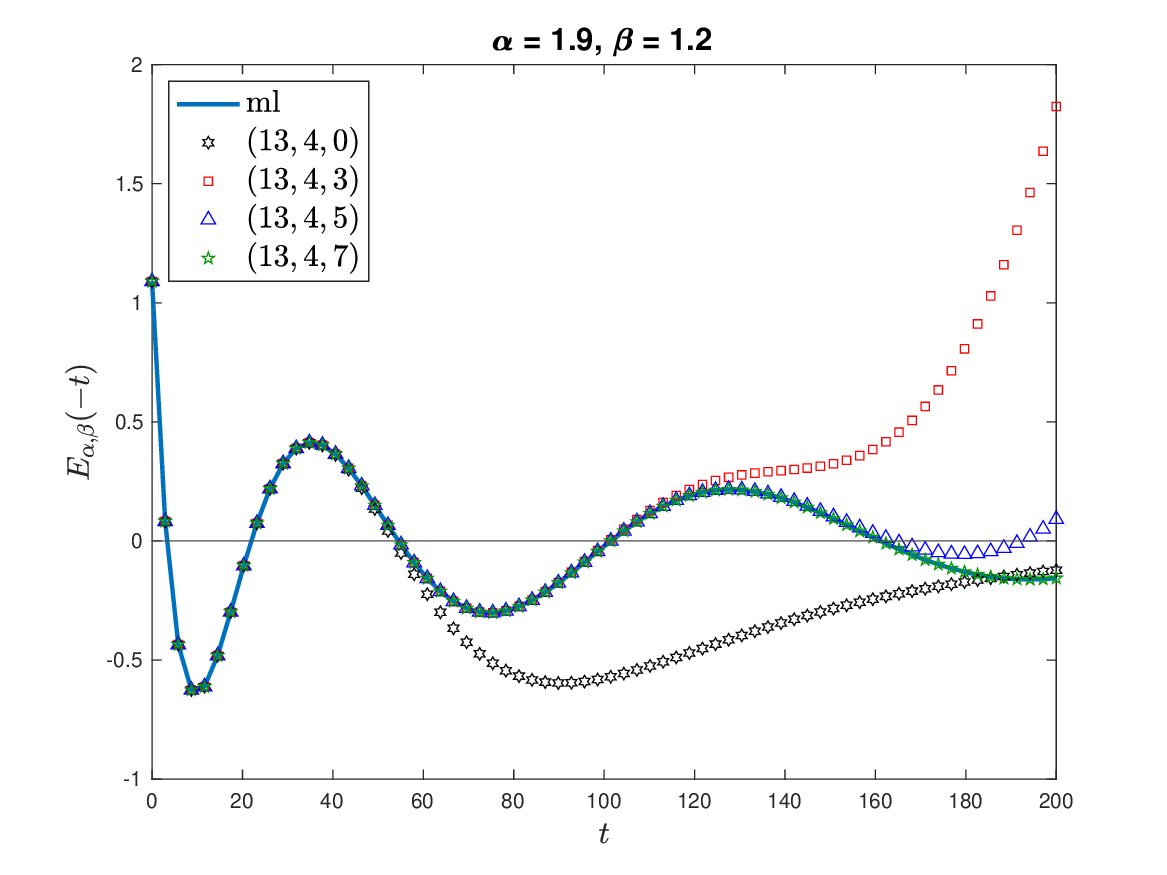}
\includegraphics[width=.49\textwidth,height=.3\textwidth]{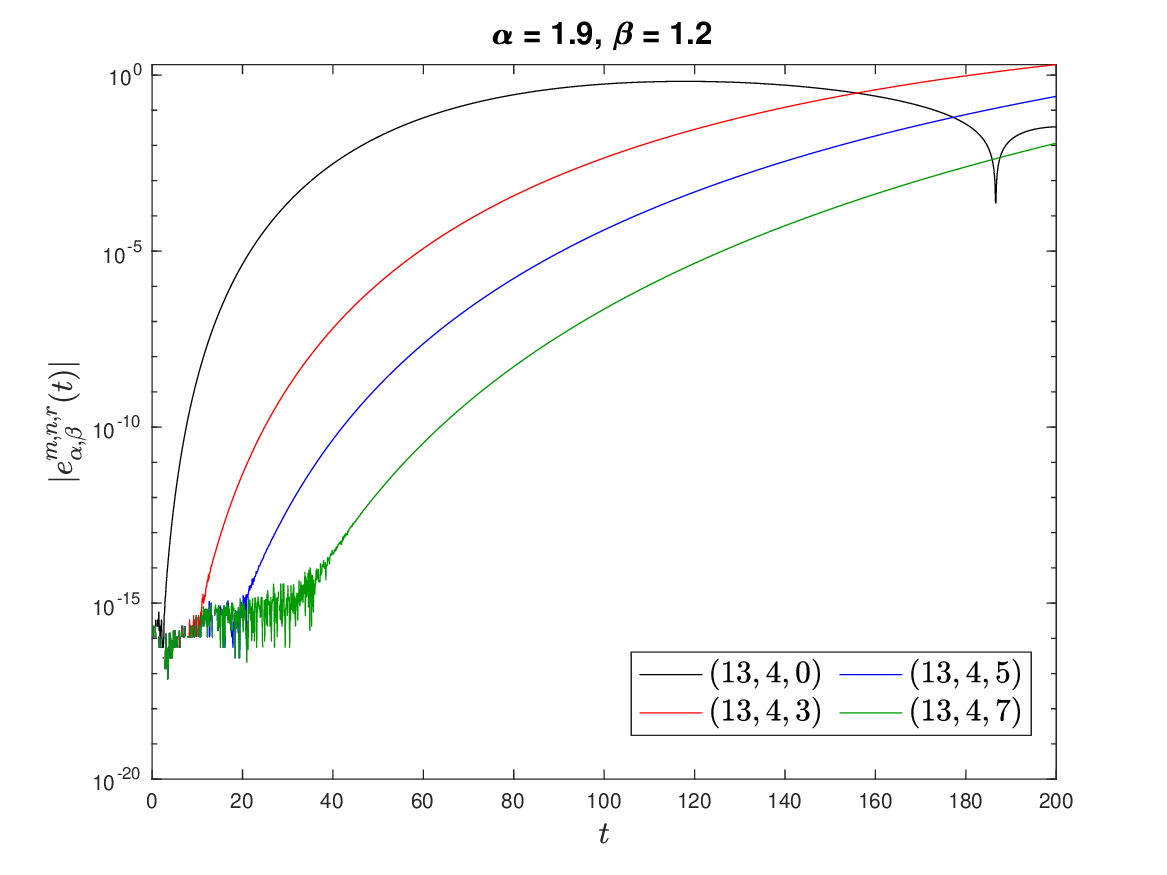}
\caption{Plots of $R_{\alpha,\beta}^{13,4, r}$ for different values of $r$
when $(\alpha,\beta)$ in region (A).}
\label{fig:7}
\end{figure}

\begin{figure}[p]
\centering
\includegraphics[width=0.49\textwidth,height=.3\textwidth]{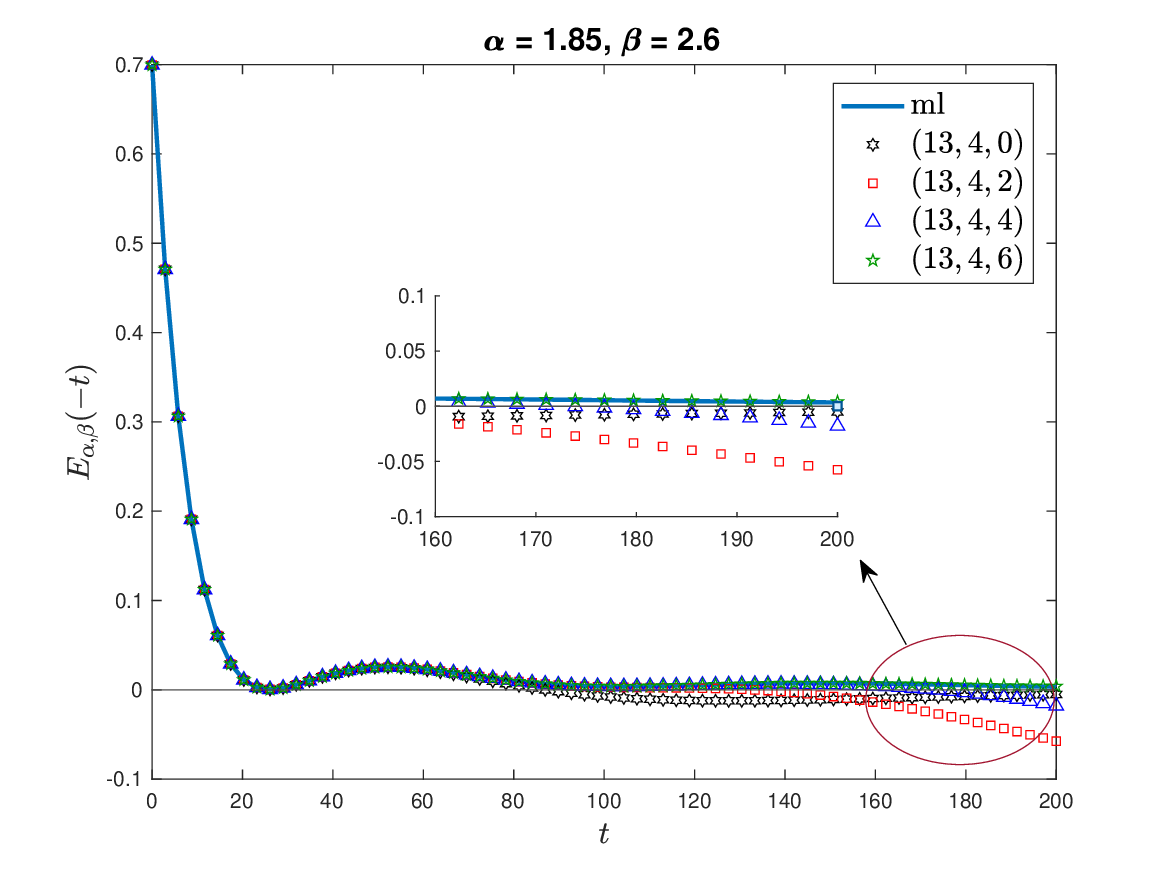}
\includegraphics[width=0.49\textwidth,height=.3\textwidth]{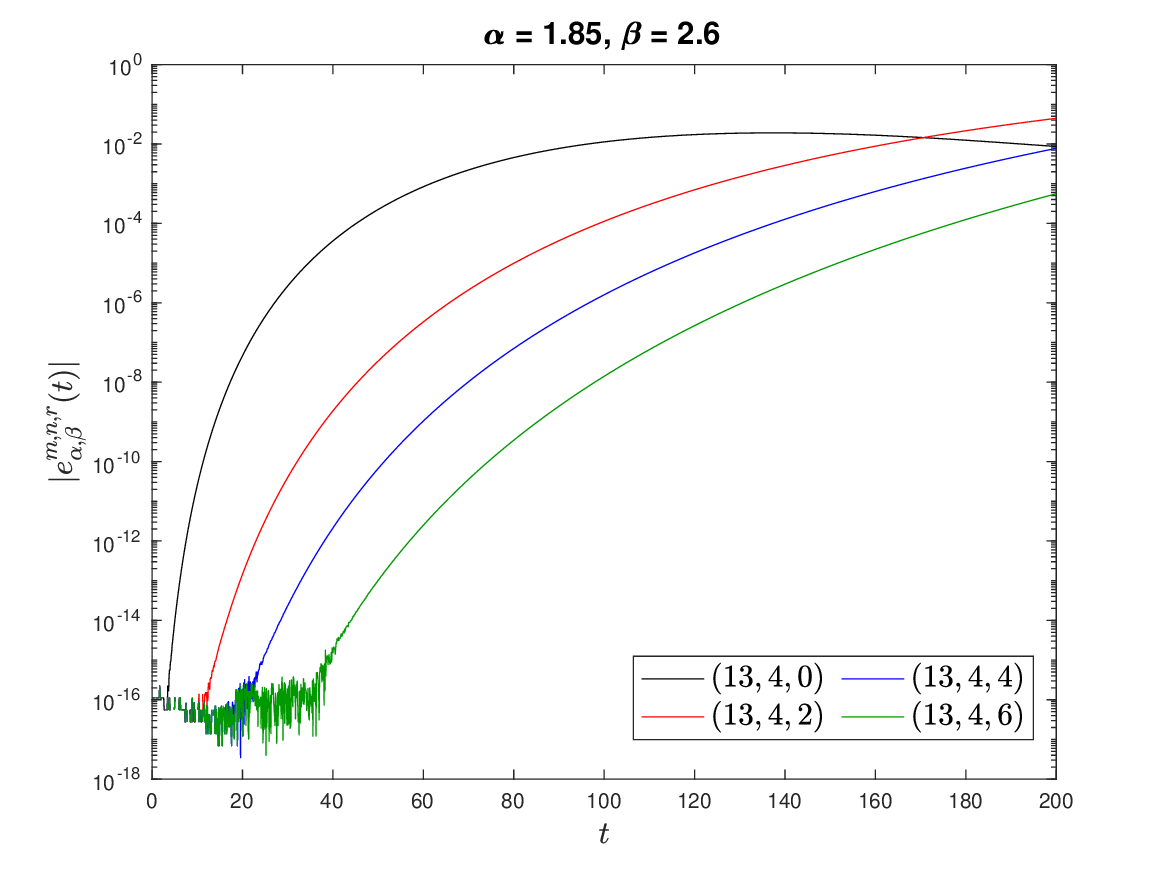}
\caption{Plots of $R_{\alpha,\beta}^{13,4, r}$ for different values of $r$
when $(\alpha,\beta)$ in region (C).}
\label{fig:6}
\end{figure}

\begin{figure}[p]
\centering
\includegraphics[width=0.49\textwidth,height=.3\textwidth]{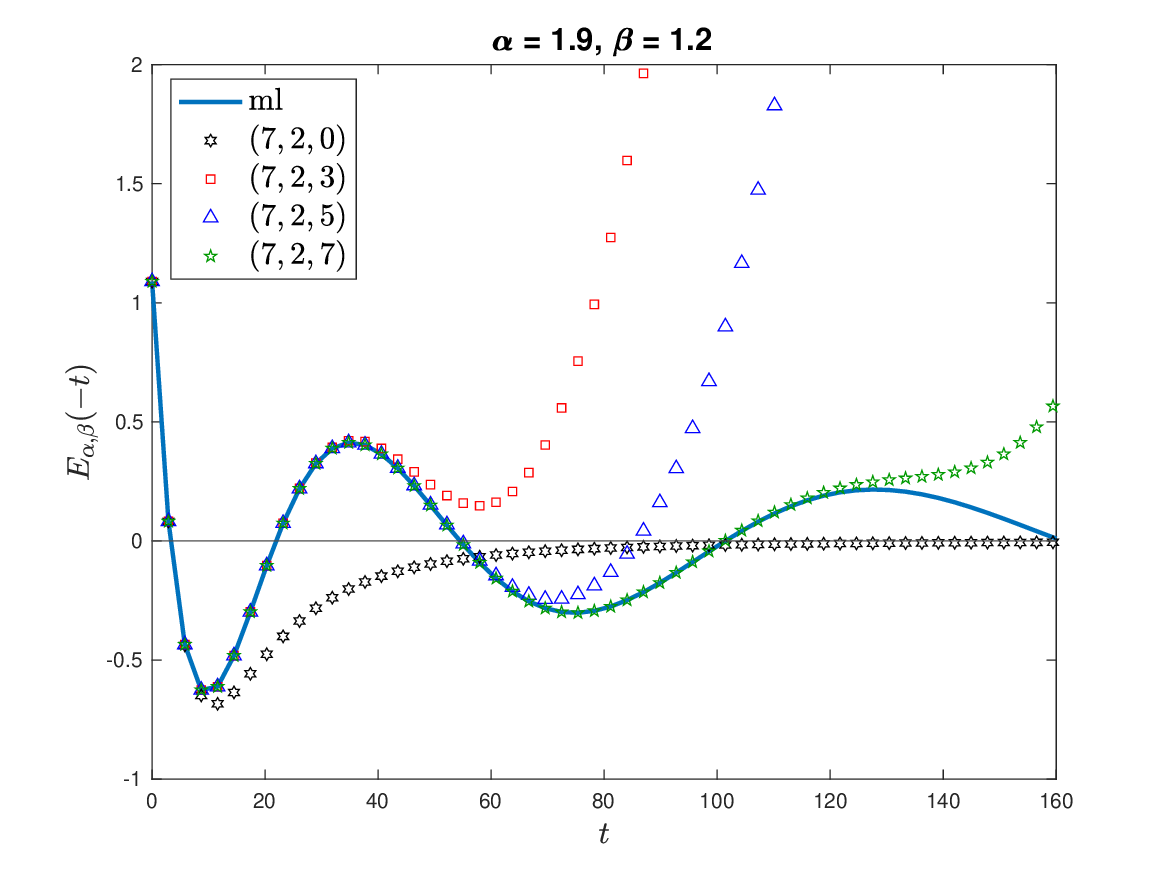}
\includegraphics[width=0.49\textwidth,height=.3\textwidth]{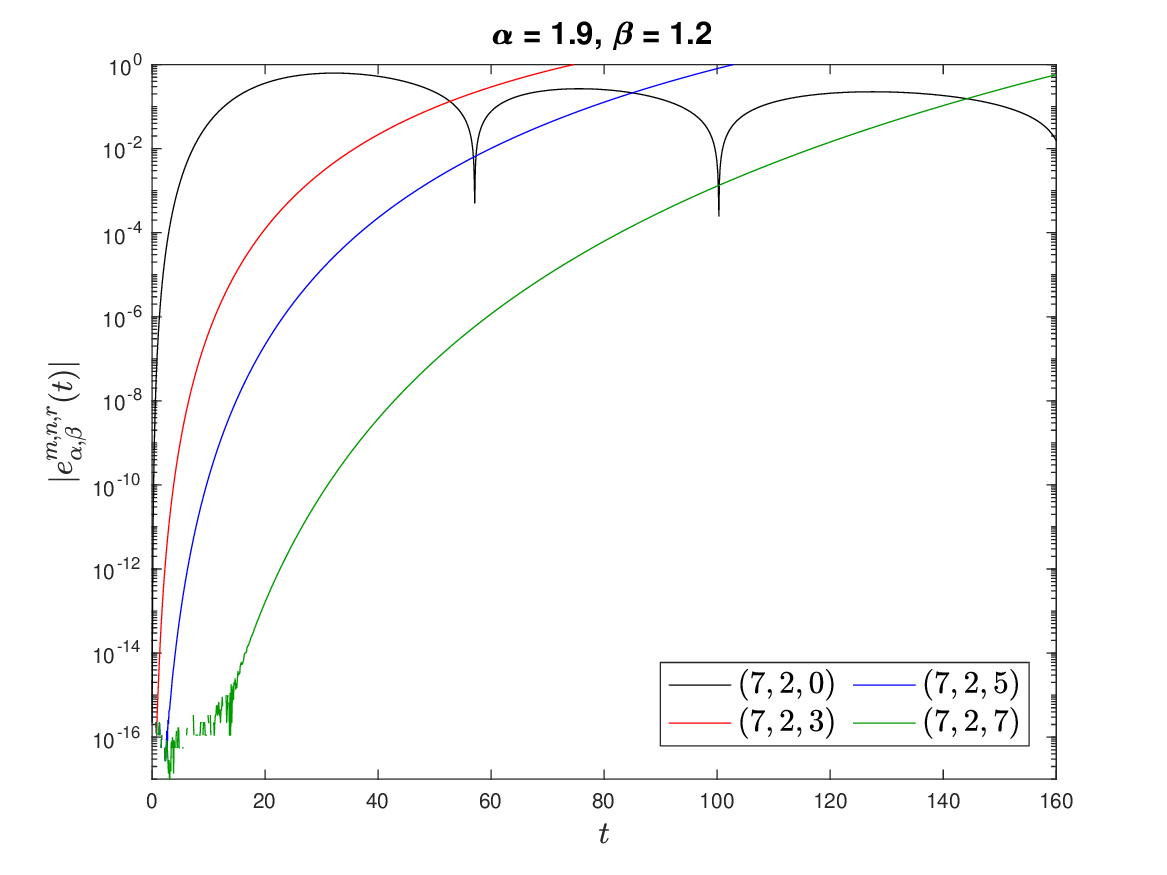}
\caption{ Plots of $R_{\alpha,\beta}^{7,2, r}$ for different values of $r$
when $(\alpha,\beta)$ in region (A).}
\label{fig:8}
\end{figure}

The performance of the generalized approximant $R_{\alpha,\beta}^{13,4,r}$ 
when $(\alpha,\beta)$ in regions (A) and (C) is demonstrated 
in Figures~\ref{fig:7} and \ref{fig:6}, respectively, where
$$
\left|e_{\alpha, \beta}^{m, n,r}(t)\right|:=\left|E_{\alpha, \beta}(-t)-R_{\alpha, \beta}^{m, n,r}(t)\right|.
$$
As can be observed, the interval of approximation can be extended by increasing $r$.
The generalized approximant $R_{\alpha,\beta}^{7,2,r}$ has the same features, although it is not as accurate as $R_{\alpha,\beta}^{14,3,r}$, as shown in Figure \ref{fig:8}.


\newpage
\section{Approximation of matrix Mittag–Leffler}
\label{sec:matrix MLF}

In this section we study the approximation of matrix MLF defined by \cite{Sadeghi2018} 
\begin{equation}
E_{\alpha, \beta}(A)=\sum_{k=0}^{\infty} \frac{A^k}{\Gamma(\alpha k+\beta)}, \quad
\operatorname{Re} \alpha>0, \beta \in \mathbb{C}, \quad A \in \mathbb{C}^{n \times n}.
\end{equation}
Using the definition of rational functions of matrices, see \cite{higham2008functions}, 
we define the global Pad\'e approximation of $E_{\alpha, \beta}(A)$ by	
\begin{equation} \label{app:def1}
R^{m,n}_{\alpha,\beta}(A) := \frac{1}{\Gamma(\beta-\alpha)} [q(A)]^{-1} p(A),
\end{equation}
for some appropriate matrix $A$, where $p$ and $q$ denote the numerator and 
denominator, respectively, of
$\Gamma(\beta-\alpha) R^{m,n}_{\alpha,\beta}$.
Furthermore, by extending the decomposition \eqref{eq:rec1} to matrix arguments, we define the 
rational approximant 
\begin{equation} \label{eq:mlf dec matrices}
R^{m,n,r}_{\alpha,\beta}(A) \defeq (-A)^r R^{m,n}_{\alpha,\beta+r\alpha}(A) 
+ P_{\alpha,\beta}^{r-1}(-A), \qquad r=0, 1, 2, \dots ,
\end{equation} 
where $P_{\alpha,\beta}^{r-1}$ is the polynomial given by 
\eqref{eq:mlf poly term}. 

Next we discuss some approaches for implementing \eqref{app:def1} and
compare their accuracy and computation time.


\begin{enumerate}

\item
Matrix Inversion approach \\
In this approach the inverse of $q(A)$ is calculated and then substituted in \eqref{app:def1}.

\item
Linear System approach \\
The approximant $R_{m,n}^{\alpha,\beta}(A)$ is obtained by solving the matrix system
$$
q(A) R_{m,n}^{\alpha,\beta}(A) = p(A).
$$
This requires solving $N$ systems for an $N\times N$ matrix.

\item 
Partial Fraction approach \\
Partial fraction decomposition is known to provide an efficient form for evaluating rational functions. 
For the global Pad\'e approximants, it has been discussed in 
\cite{Sarumi2020,Sarumi2021} that these approximants have complex conjugate roots which can add to the efficient implementation. 
As an example, the approximant 
$R_{\alpha, \beta}^{13,4}(x)$ which admits the partial fraction decomposition
\begin{equation}
R_{\alpha, \beta}^{13,4}(x) = \sum_{i=1}^4 
\left[ \frac{c_i}{x-s_i} + \frac{\bar{c}_i}{x-\bar{s}_i} \right],
\end{equation}
where $\left\{c_1, c_2, c_3, c_4\right\}$ and $\left\{s_1, s_2, s_3, s_4\right\}$ are the non-conjugate residues and poles, respectively, can be written as
\begin{equation}
R_{\alpha, \beta}^{13,4}(x) = 2 \operatorname{Re} \sum_{i=1}^4 
\left[ \frac{c_i}{x-s_i} \right].
\end{equation}
So for a matrix argument $A$, the approximant can be calculated as
$$
E_{\alpha, \beta}(-A)  \approx R_{\alpha, \beta}^{13,4}(A)  
\defeq 2 \operatorname{Re} \sum_{j=1}^4 c_j\left(A-s_j \mathrm{I}\right)^{-1},
$$
where $I$ is the identity matrix.  

\item
Matrix Diagonalization approach \\
When the matrix argument $A$ is diagonalizable, then the factorization 
$A = Z D Z^{-1}$ could be considered, where $D$ is the diagonal matrix containing the eigenvalues and the columns of $Z$ are the corresponding eigenvectors. In this case, the Matrix MLF can be computed as 
$$
E_{\alpha,\beta}(A) = Z\, E_{\alpha,\beta}(D)\, Z^{-1} = Z \operatorname{diag}(E_{\alpha,\beta}(\lambda_i))\, Z^{-1}.
$$
Accordingly, the approximant $R^{m,n}_{\alpha,\beta}$ could be computed as 
$$
R^{m,n}_{\alpha,\beta}(A) = 
Z \operatorname{diag}(R^{m,n}_{\alpha,\beta}(\lambda_i))\, Z^{-1} .
$$
\end{enumerate}

\begin{figure}
\centering
\includegraphics[width=0.49\textwidth]{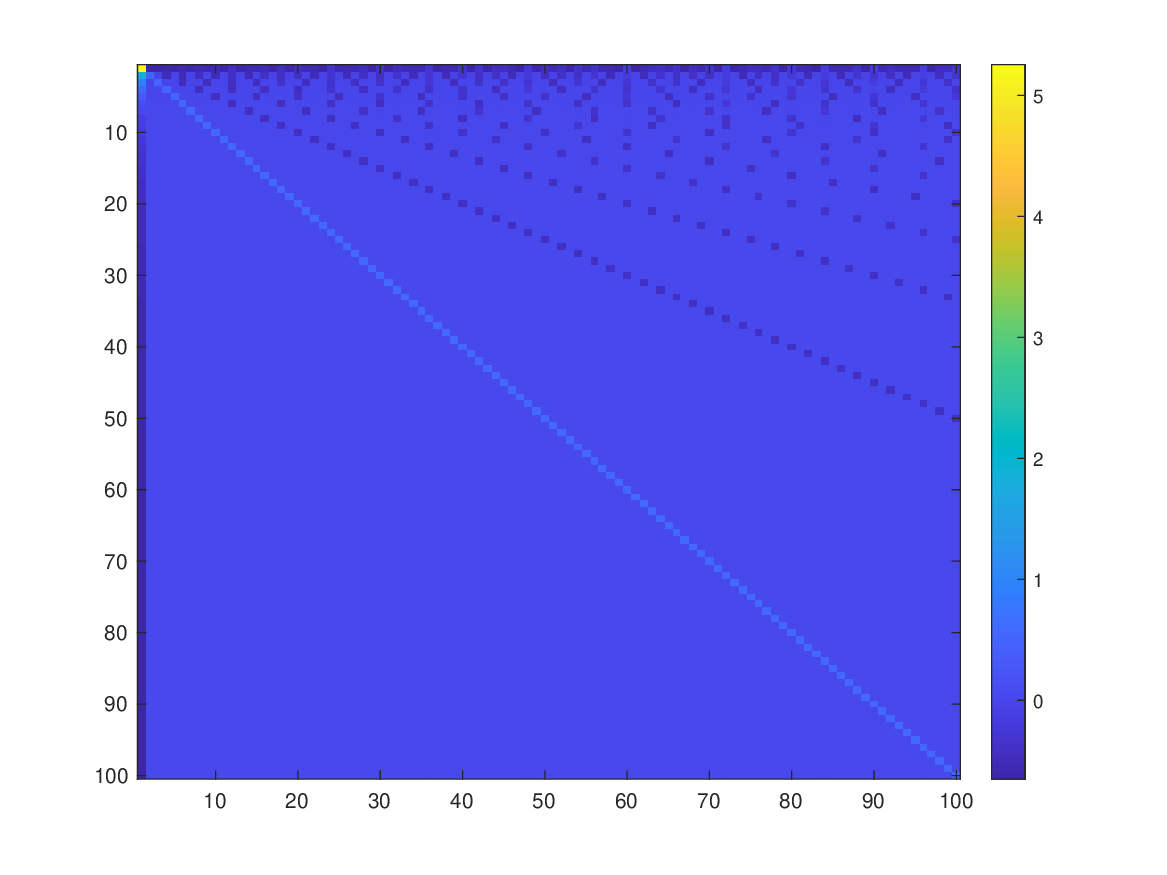} 
\includegraphics[width=0.49\textwidth]{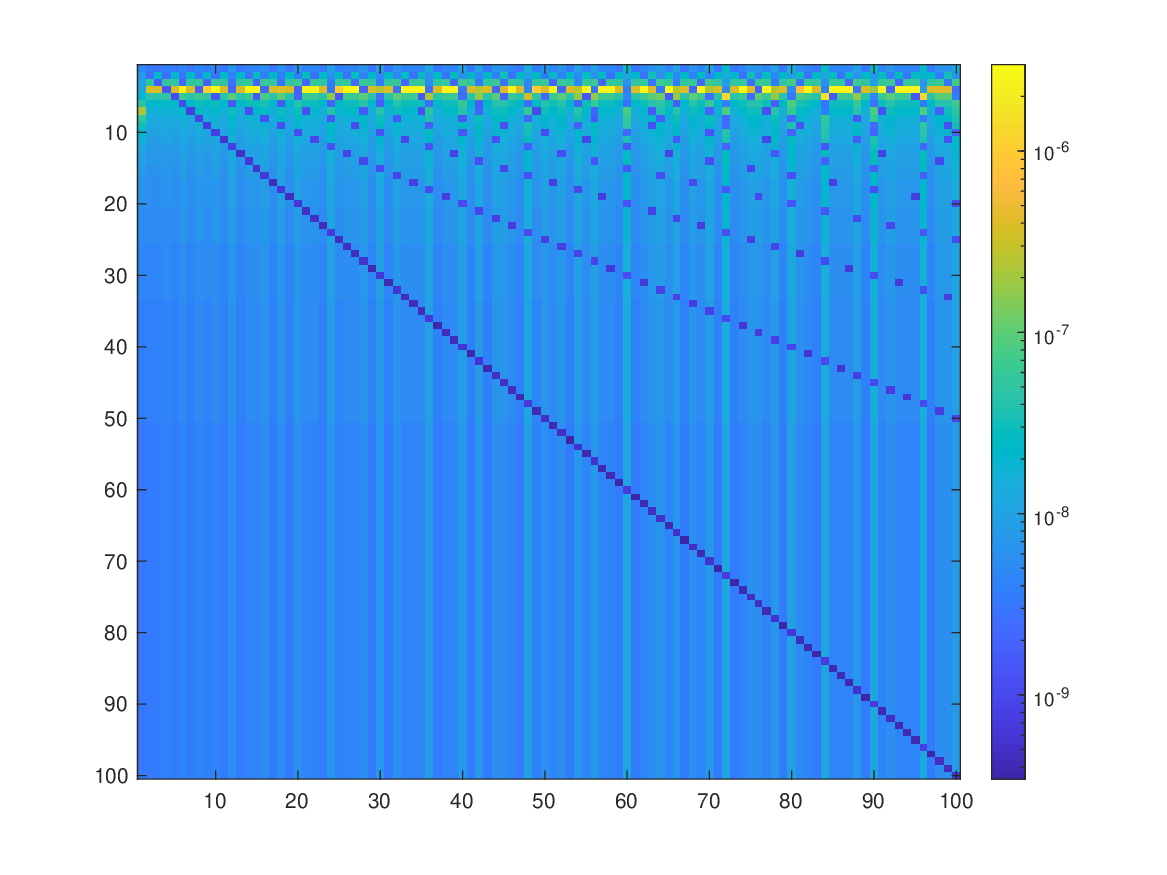}
\caption{MLF matrix approximation $R^{13,4}_{\alpha,\beta}$, 
$(\alpha,\beta)= (1.9, 1)$, of 100$\times$100 Redheffer matrix using 
partial fractions (Left) and its relative error (Right).}
\label{fig:redheff1}
\end{figure}

For experimental purposes, we consider the performance of the approximant 
\eqref{eq:general} when applied to Redheffer matrix of size 100$\times$100,
in comparison to the reference values by ml\_matrix. 
In Figure \ref{fig:redheff1}, a color map is showing the component-wise values and relative errors for $R^{13,4}_{\alpha,\beta}$, $(\alpha,\beta)=(1.9, 1)$,
when calculated using partial fractions.
In Table \ref{tab:redheff}, a comparison of the absolute error, relative error, and runtime is provided for the different techniques in computing $R^{13,4}_{\alpha,\beta}$.

\renewcommand{\arraystretch}{1.5}
\begin{table}[]
\centering
\caption {Maximum absolute error (AE), maximum relative error (RE) and the runtime (in seconds) for computing the matrix MLF of a $100 \times 100$ Redheffer matrix, $\beta = 1$.} 
\label{tab:redheff}  
\begin{tabular}{l|ccc|ccc}
& \multicolumn{3}{c}{$\alpha = 1.5$} & \multicolumn{3}{c}{$\alpha = 1.9$} 
\\
& AE         & RE        & Runtime   & AE         & RE        & Runtime  
\\  \hline
Matrix Inversion & 1.13E-04   & 7.20E-03  & 2.00E-03  & 2.41E-08   & 3.00E-06  & 1.10E-03  \\
Linear System    & 1.13E-04   & 7.20E-03  & 1.90E-03  & 2.41E-08   & 3.00E-06  & 1.10E-03  \\
Partial Fraction & 1.13E-04   & 7.20E-03  & 3.70E-03  & 2.41E-08   & 3.00E-06  & 2.30E-03  \\
ml\_matrix       &            &           & 1.44E-01  &            &           & 9.73E+00 
\end{tabular}
\end{table}

\section{Applications and numerical experiments}
\label{sec:applications}

We present here some applications of the general Pad\'e approximant \eqref{eq:general} related to fractional oscillation equations where MLFs arise naturally as solutions. 
Numerical experiments are included to highlight the efficiency and accuracy of these approximants.

\subsection{\bf Application: Fractional plasma oscillations}
\label{exp:fracPlas}

The fractional plasma oscillation model, as in \cite{aguilar2014fractional}, is given by 
\begin{equation} \label{eq:model}
\begin{gathered}
\cD^{\alpha} u(t)+ Au(t) = f(t),\\  
u(0)=u_{0},\\
u^{\prime}(0)=u_{1},
\end{gathered}
\end{equation}
where $\cD^{\alpha}$ denotes the Caputo fractional derivative of order $\alpha \in(1,2)$ defined by
\[
\cD^{\alpha}u(t):= \begin{cases}
\dfrac{1}{\Gamma(2-\alpha)} \displaystyle \int_{0}^{t} {(t-\tau)^{1-\alpha}} u^{\prime\prime}(\tau) d \tau, & 1<\alpha<2,\\
u^{\prime\prime}(t), & \alpha=2.
\end{cases}
\]
The constant $A$ is the fractional electron plasma frequency and $f(t)$ is the electric field.
We consider the model with static electric field and with no electric field.

\subsubsection
{\bf Fractional plasma oscillations model with static electric field}
\label{subsub:static}

Consider the special case of problem \eqref{eq:model}, 
$$
\cD^{\alpha} u(t)+ u(t) = 1, \qquad 
u(0) = 1, \qquad 	u^{\prime}(0) = -1.
$$
The exact solution for this initial value problem is   
\begin{equation}
\label{eq:fracplas solution}
u(t)=E_{\alpha}\left(-t^{\alpha}\right) 
- t E_{\alpha, 2}\left(-t^{\alpha}\right) 
+ t^{\alpha} E_{\alpha, \alpha+1}\left(-t^{\alpha}\right).
\end{equation}

\begin{figure}[p]
\centering
\includegraphics[width=0.49\textwidth]{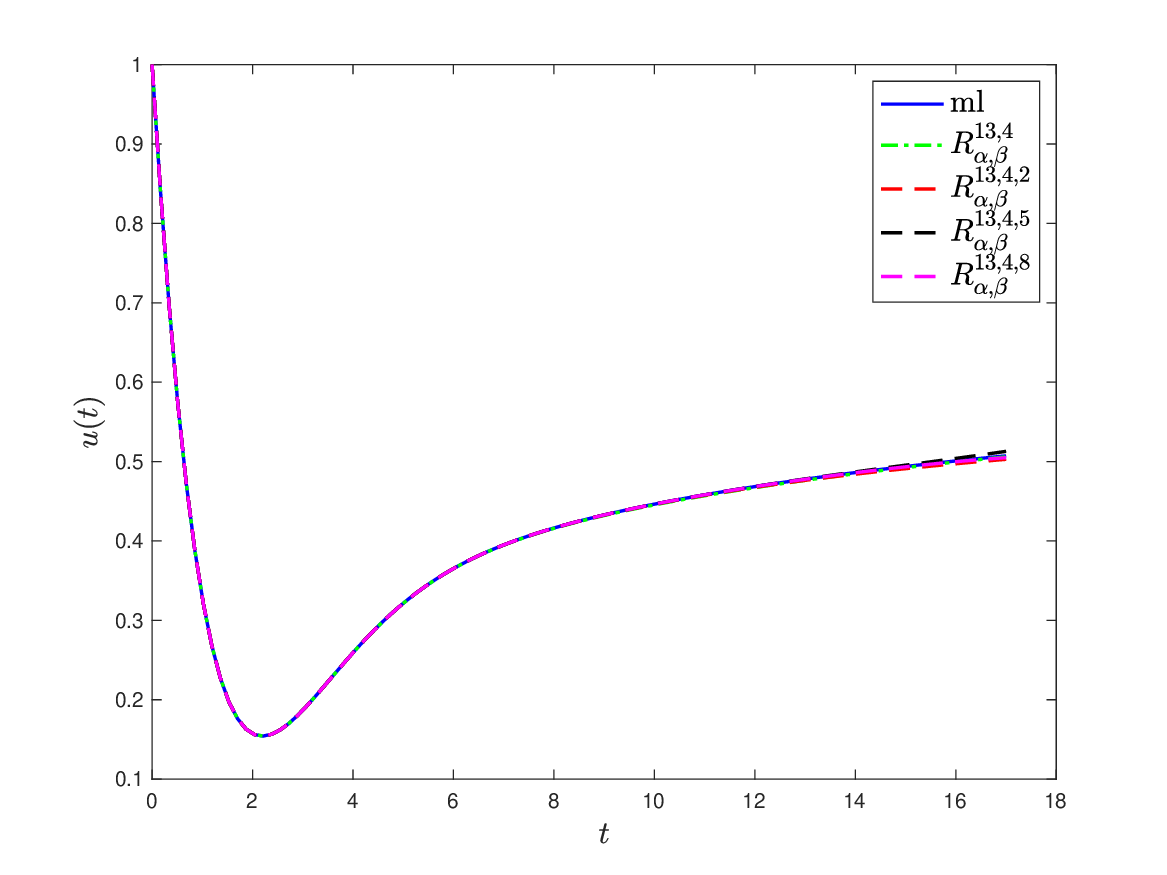}
\includegraphics[width=0.49\textwidth]{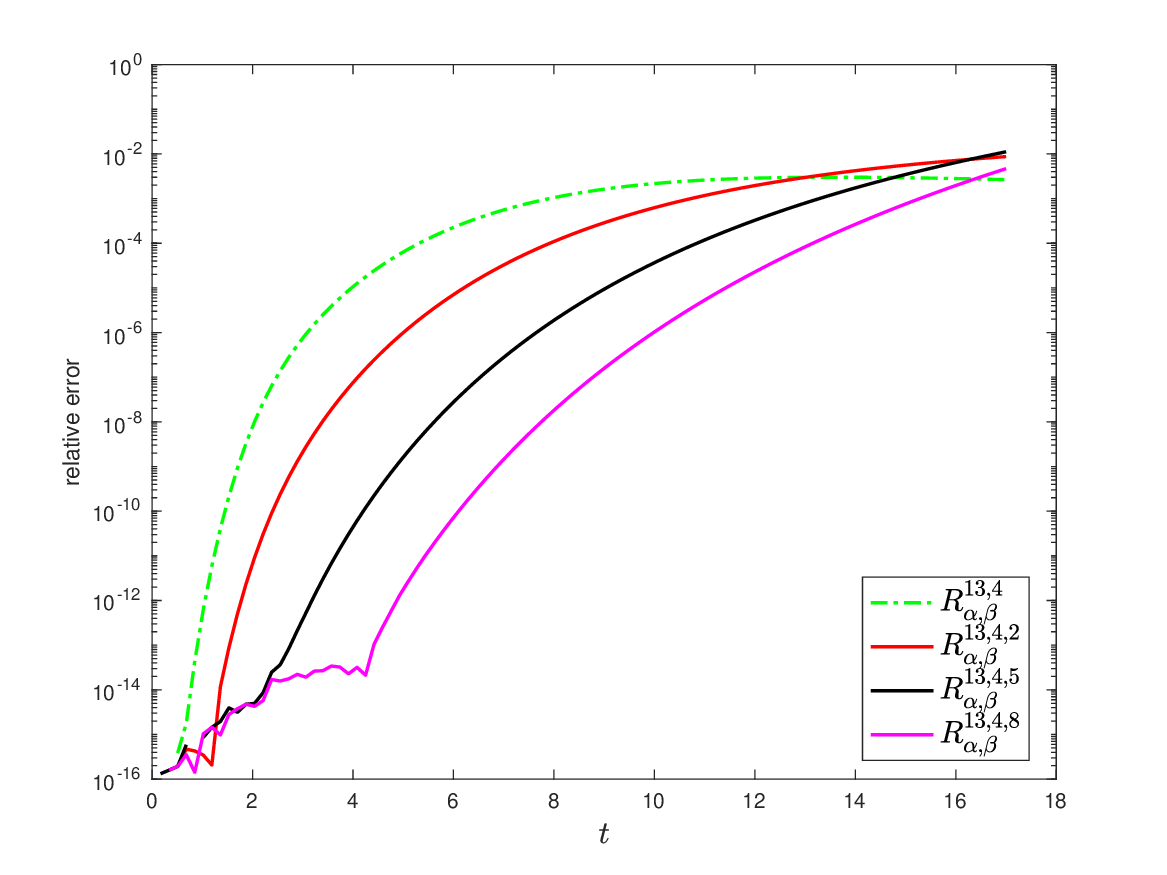}
\caption{Profiles (left) and relative error (right) in approximating the solution
\eqref{eq:fracplas solution}, $\alpha=1.2$.}
\label{fig:fracPlas_a12}
\end{figure}

\begin{figure}[p]
\centering
\includegraphics[width=0.49\textwidth]{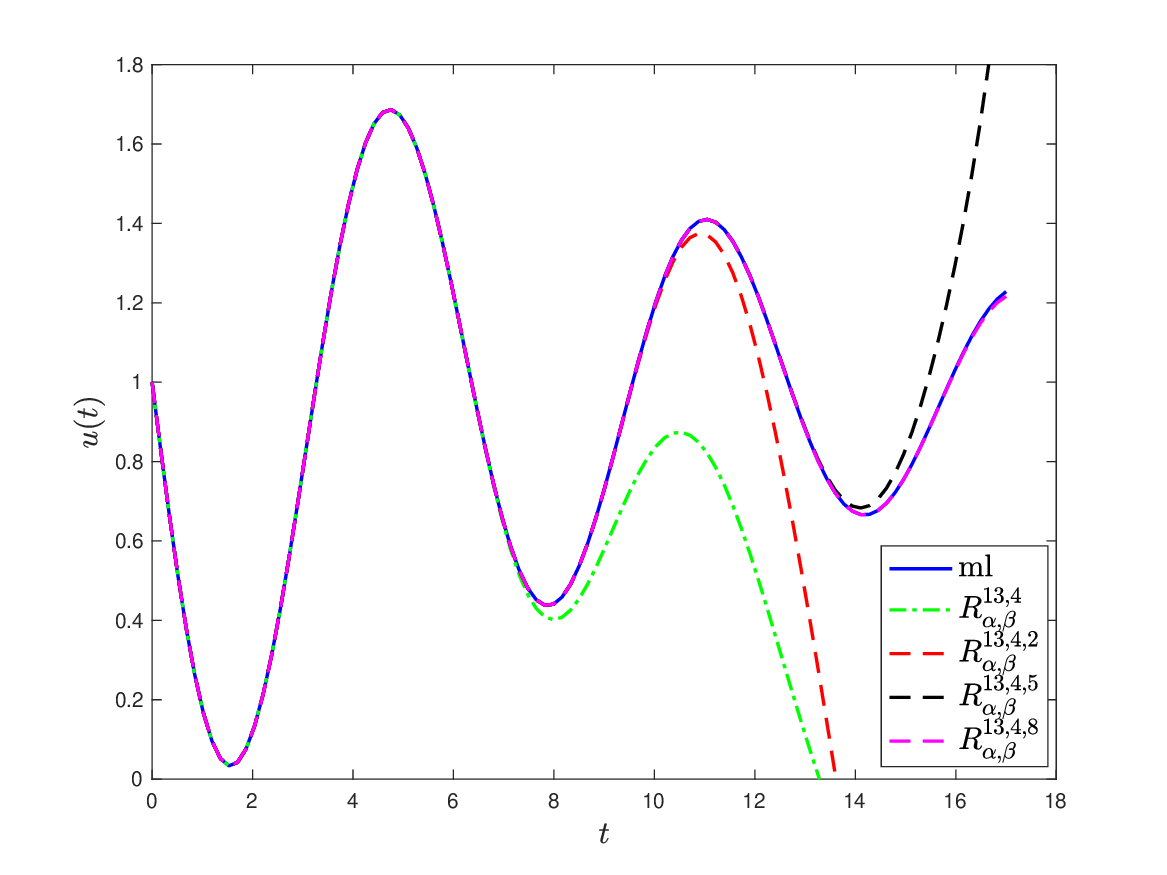}
\includegraphics[width=0.49\textwidth]{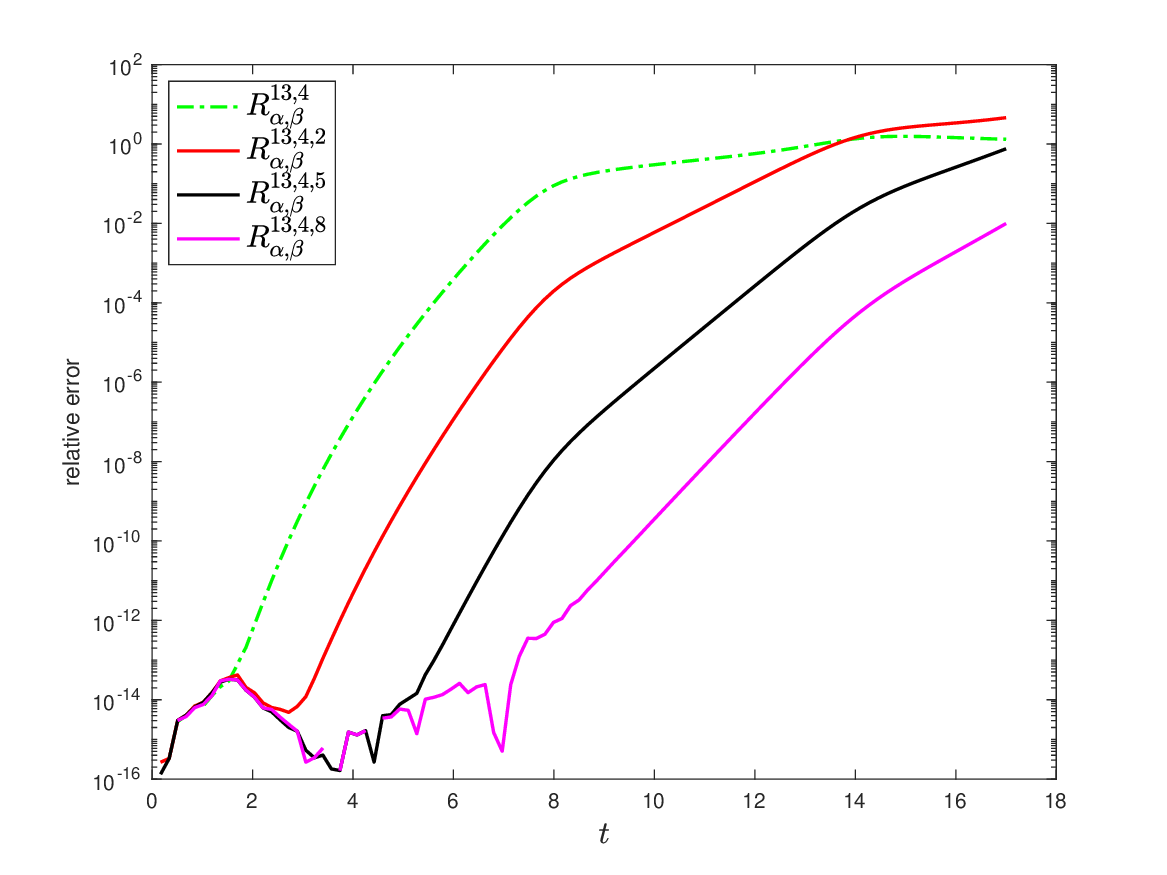}
\caption{Profiles (left) and relative error (right) in approximating the solution
\eqref{eq:fracplas solution}, $\alpha=1.9$.}
\label{fig:fracPlas_a19}
\end{figure}

In Figure \ref{fig:fracPlas_a12}, a comparison of approximations of \eqref{eq:fracplas solution} when $\alpha = 1.2$ is provided.
Considering the three terms in \eqref{eq:fracplas solution},
the pair $(\alpha,\beta)$ when $\beta=1$ falls in region (A)
of the phase diagram \ref{fig:diag}, where the MLF is oscillatory. 
While the pairs $(1.2, 2)$ and $(1.2, 2.2)$ fall 
in region (D) and (F), respectively, in which the MLF has 
no real roots and no oscillations.
A similar comparison when $\alpha = 1.9$ is included in 
Figure \ref{fig:fracPlas_a19}.
In this case, the terms $E_{1.9,1}$ and $E_{1.9,2}$ correspond to $(\alpha,\beta)$ in parts of Region (A) in which the MLF is highly oscillatory and has many real roots. 

To illustrate the computation efficiency, it is seen in 
Table \ref{tab:fracPlas} that $R_{\alpha,\beta}^{13,4}$ and  $R_{\alpha,\beta}^{13,4,r}$ with $r=2,5,8$ have 
significantly less computation cost (in time) than the ml function. Moreover, while the 
modified approximant $R_{\alpha,\beta}^{13,4,8}$ yields more accurate approximations than 
$R_{\alpha,\beta}^{13,4}$ when compared to the reference values, they have comparable 
computation time.

\begin{table}[]
	\centering 
\caption {Maximum relative error (RE) and runtime in computing the solution 
	\eqref{eq:fracplas solution}, over the time mesh $[0, 0.17,0.34, \dots, 17]$,  $\alpha = 1.9$. }
\label{tab:fracPlas}
\begin{tabular}{lcc}
	& RE       & Runtime  \\ \hline
	$R^{13,4}_{\alpha,\beta}$   & 1.56E+00 &5.69E-04 \\
	$R^{13,4,2}_{\alpha,\beta}$ & 4.59E+00 & 5.19E-04 \\
	$R^{13,4,5}_{\alpha,\beta}$ &7.56E-01 & 5.20E-04\\
	$R^{13,4,8}_{\alpha,\beta}$ & 9.91E-03 & 5.52E-04 \\
	ml            & -        & 2.38E-02
\end{tabular}
\end{table}

\subsubsection
{\bf Fractional plasma oscillations model with no electric field}

Consider the initial values problem
$$
\cD^{\alpha} u(t)+ u(t) = 0, \qquad 
u(0) = 0.2, \qquad 	u^{\prime}(0) = 0.1.
$$
The exact solution of this problem is given by 
\begin{equation}
\label{eq:fracplas-nofield}
u(t)= 0.2 E_{\alpha}\left(-t^{\alpha}\right)
+ 0.1 t E_{\alpha, 2}\left(-t^{\alpha}\right).
\end{equation}

\begin{figure}
\centering
\includegraphics[width=0.49\textwidth]{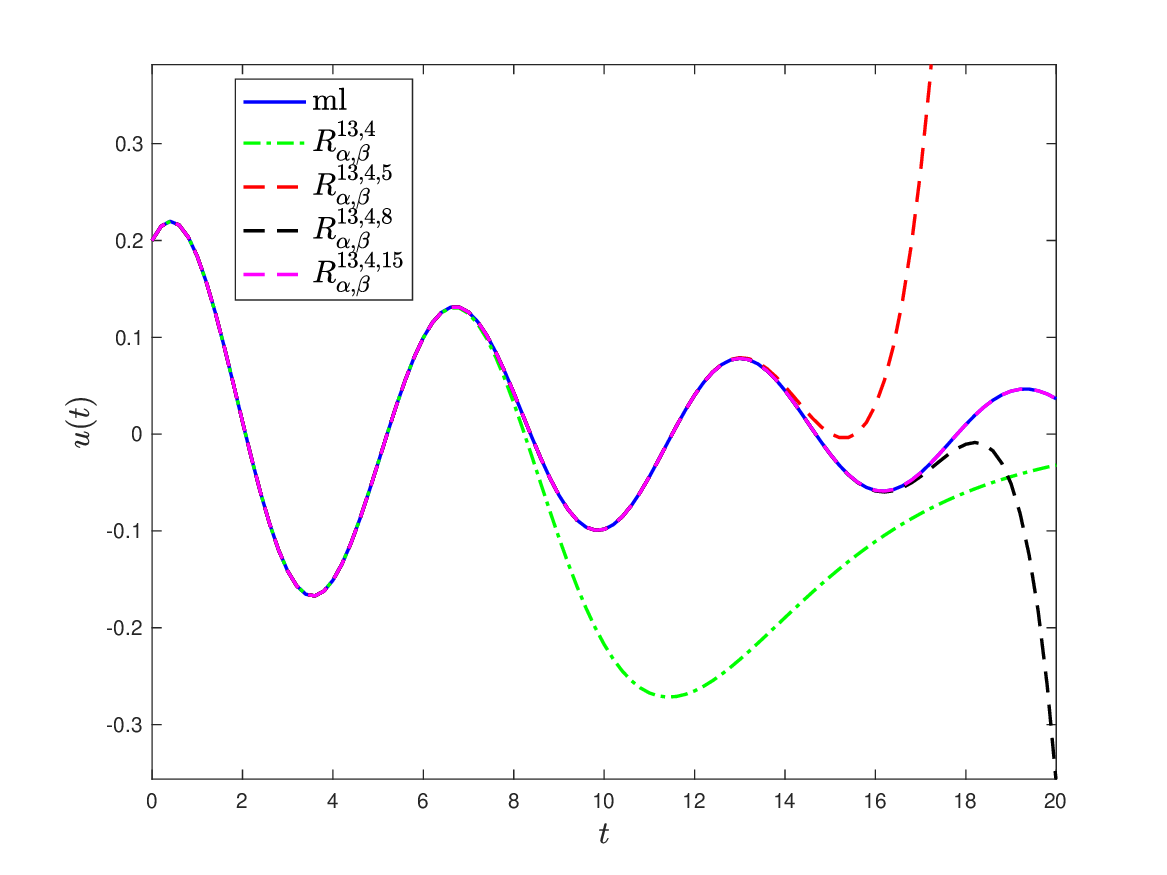}
\includegraphics[width=0.49\textwidth]{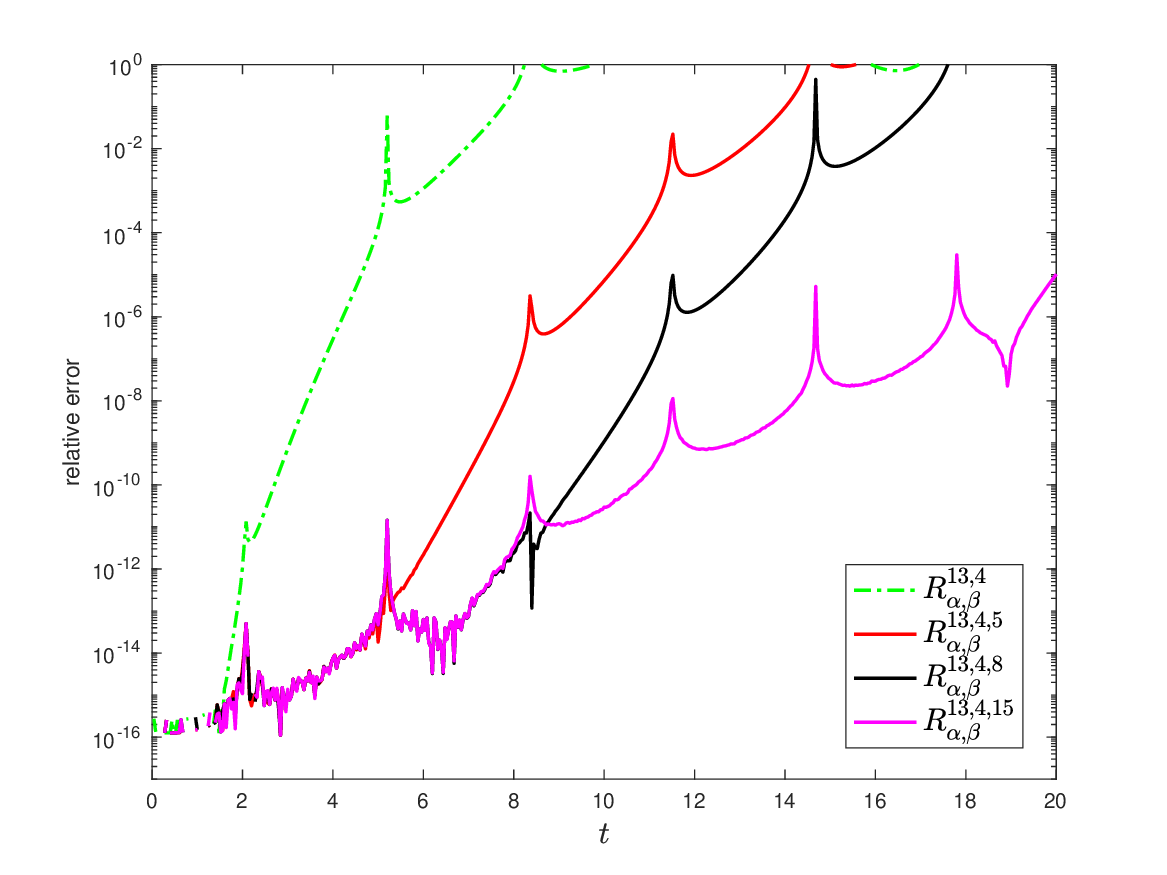}
\caption{Profiles (left) and relative error (right) in approximating the solution
\eqref{eq:fracplas-nofield}, $\alpha=1.9$.}
\label{fig:fracOsc_a19}
\end{figure}

A comparison of approximations of \eqref{eq:fracplas-nofield} when $\alpha = 1.9$ 
is provided in Figure \ref{fig:fracOsc_a19}. 
Unlike the static electric field case, it is observed that 
over an extended time interval, the solution of the fractional oscillation equation has more oscillations around zero. To sufficiently 
capture these oscillations on the interval $[0,20]$, the approximant $R^{13,4,15}_{\alpha,\beta}$ can be used.

\subsection{\bf Application: Time fractional diffusion-wave equation}
\label{wave:ex}

Consider the fractional diffusion-wave problem, \cite{murillo2011explicit}:
\begin{equation}\label{wave:ex1}
\begin{aligned}
& \cd_t^\alpha u(x, t)=\partial_x^2 u(x, t), \quad0 \leq x \leq \pi, \\
& u(x, 0)=\sin (x),  \\
&\partial_t u(x, t)|_{t=0}=0,\\
& u(0, t)=u(\pi, t)=0, \quad t \in(0, T],
\end{aligned}
\end{equation}
where $\cd_t^\alpha$ and $\partial_t$ are the fractional partial time 
derivative of order $\alpha \in(1,2)$ and 
first partial derivative in $t$, respectively, 
and $\partial_x^2$ is the second partial derivative in space. 

The exact solution of \eqref{wave:ex1} is $u(x, t) = E_\alpha(-t^{\alpha}) 
\sin (x)$, however, this problem is intended to illustrate the performance of 
the approximants \eqref{eq:general}. 
To this end, applying the second central difference approximation to (\ref{wave:ex1}) over uniform 
spatial mesh with step size $h=\pi /(m+1)$, leads to a system of the form 
\begin{equation} 
\label{eq:wave system}
\begin{aligned}
	\cD^\alpha U(t) + A U(t) & =0, \quad t>0, \\
	U(0) & =U_0,\\
	U^\prime(0) & = 0,
\end{aligned}
\end{equation}
where the $m \times m$ matrix $A=\left[a_{i j}\right]$ is tridiagonal with $a_{i, i}=2$ and $a_{i, j}=-1$ 
for $j=i \pm 1$. Further, 
$U(t) = [u_1(t), u_2(t), \dots u_m(t)]^T$, $U(0) = [\sin(x_1), \sin(x_2), \dots \sin(x_{m})]^T$, 
where $u_1(t) = u(x_1,t)$ and $x_i = i\times\pi/m$. 
The solution of the above system of differential equations is given by $U(t) = E_\alpha(-At^\alpha) U_0$.
As an experiment, we take $h = \pi/100$, so the coefficient matrix $A$ is of size $99 \times 99$. 
In addition, the rational approximation of the matrix MLF is computed by the diagonalization approach.  

\begin{figure}
\centering
\includegraphics[width=0.49\textwidth]{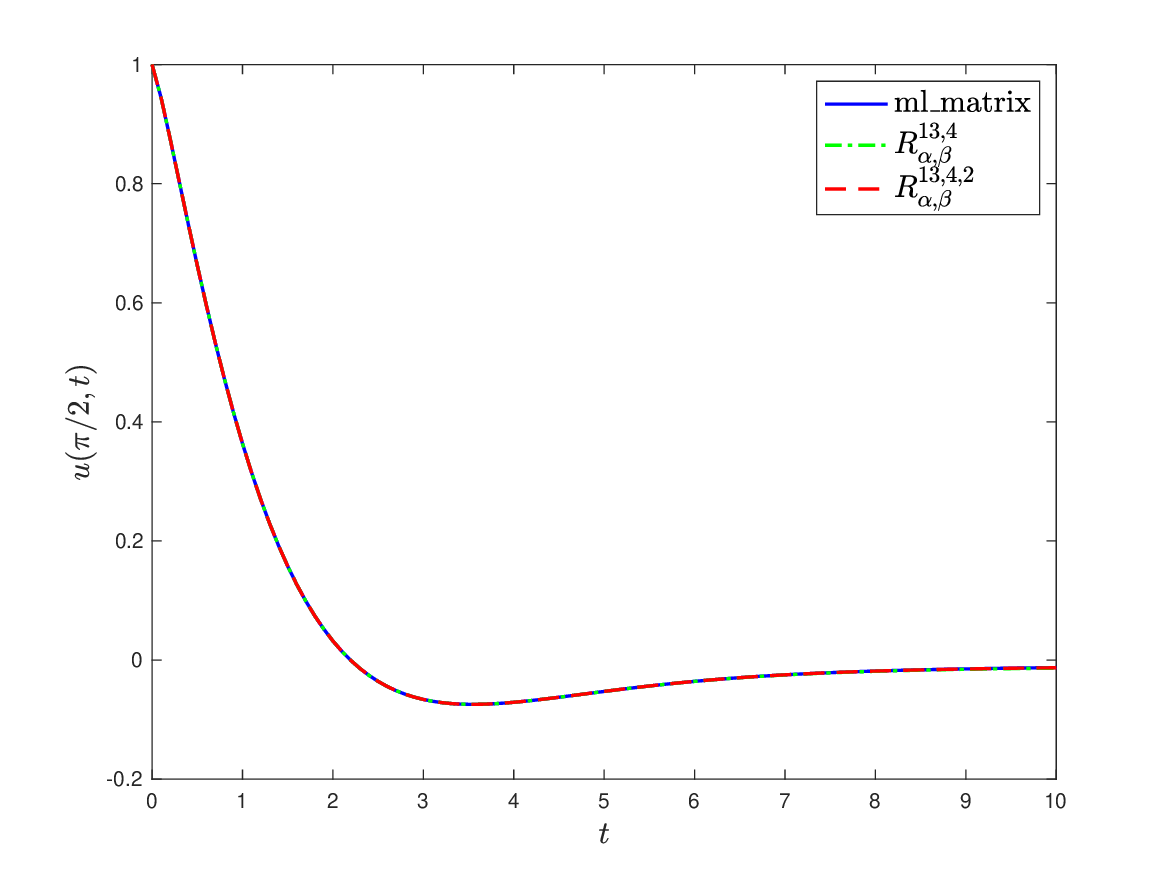}
\includegraphics[width=0.49\textwidth]{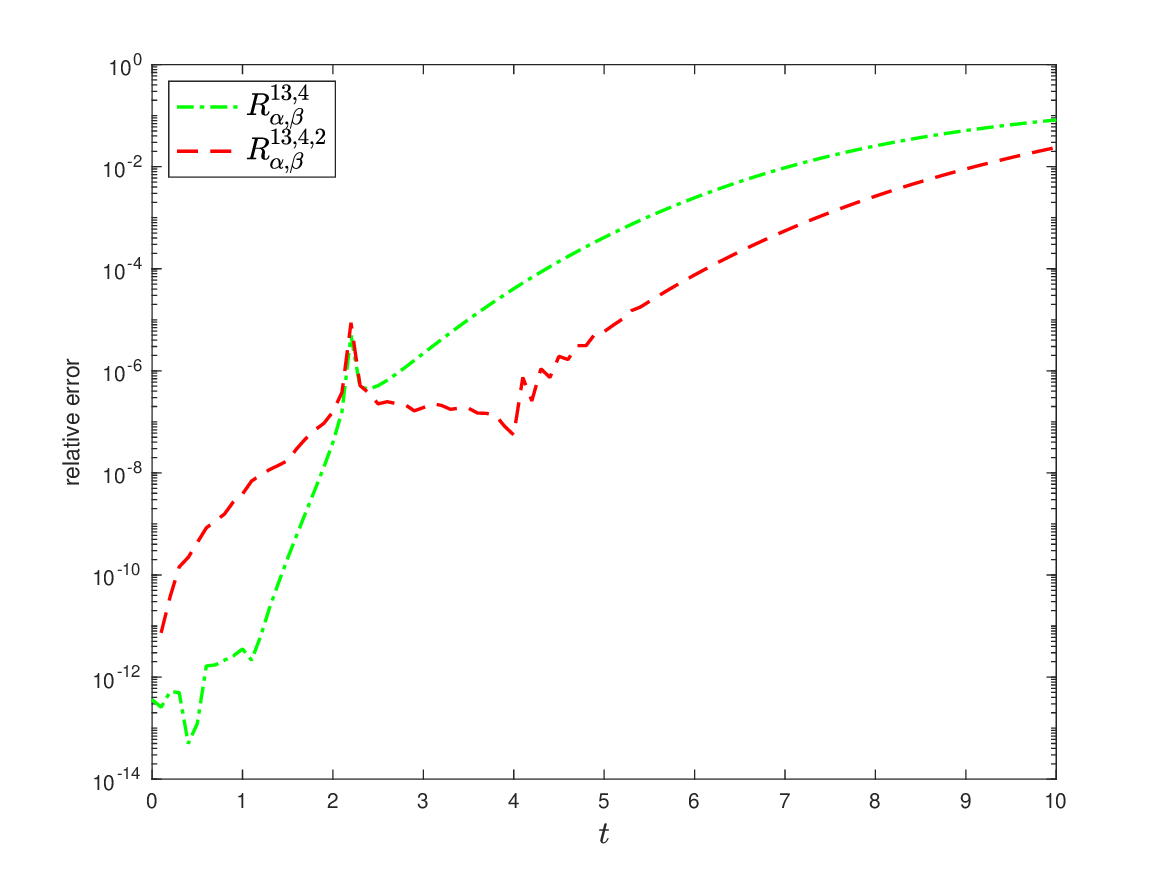}
\caption{Application \ref{wave:ex}: 
approximation of $u(\pi/2,t)$ (left) and relative error 
(right), $\alpha=1.2$.}
\label{fig:wave_a12}
\end{figure}

\begin{figure}
\centering
\includegraphics[width=0.49\textwidth]{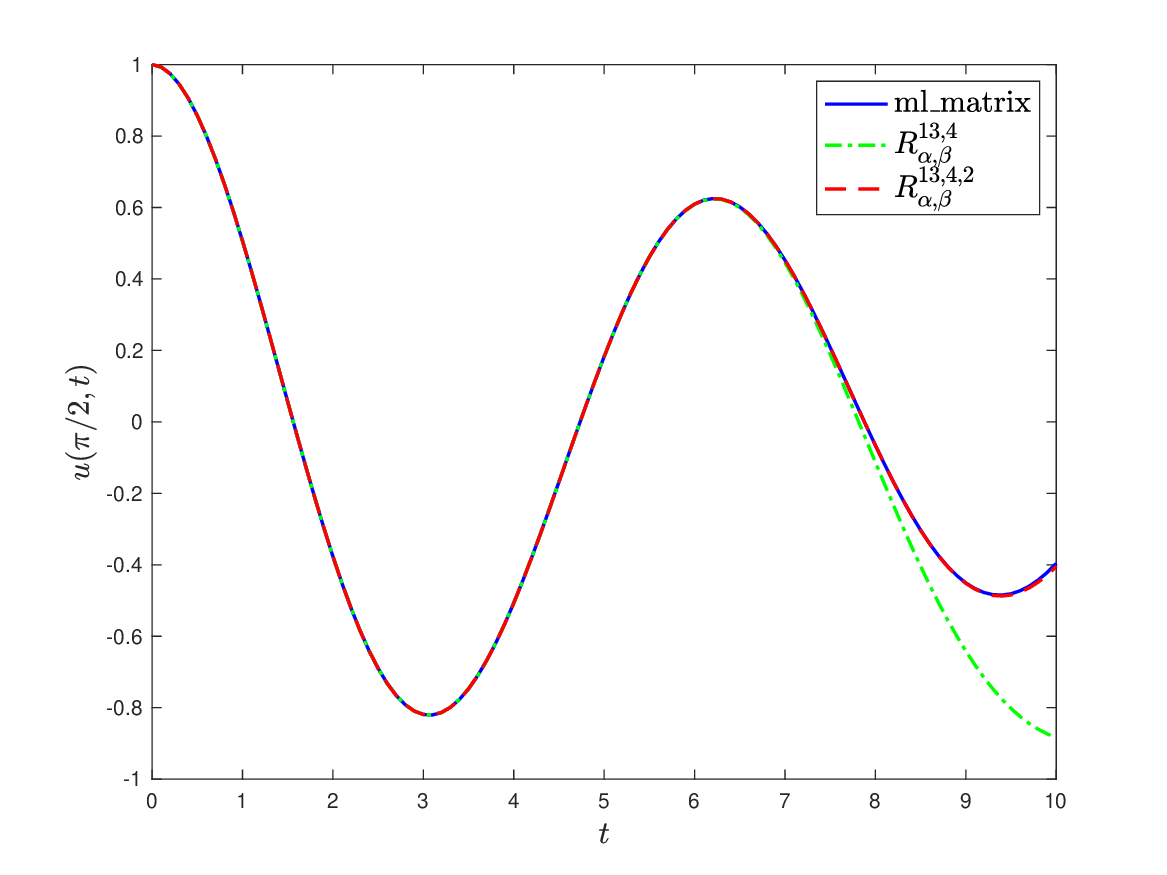}
\includegraphics[width=0.49\textwidth]{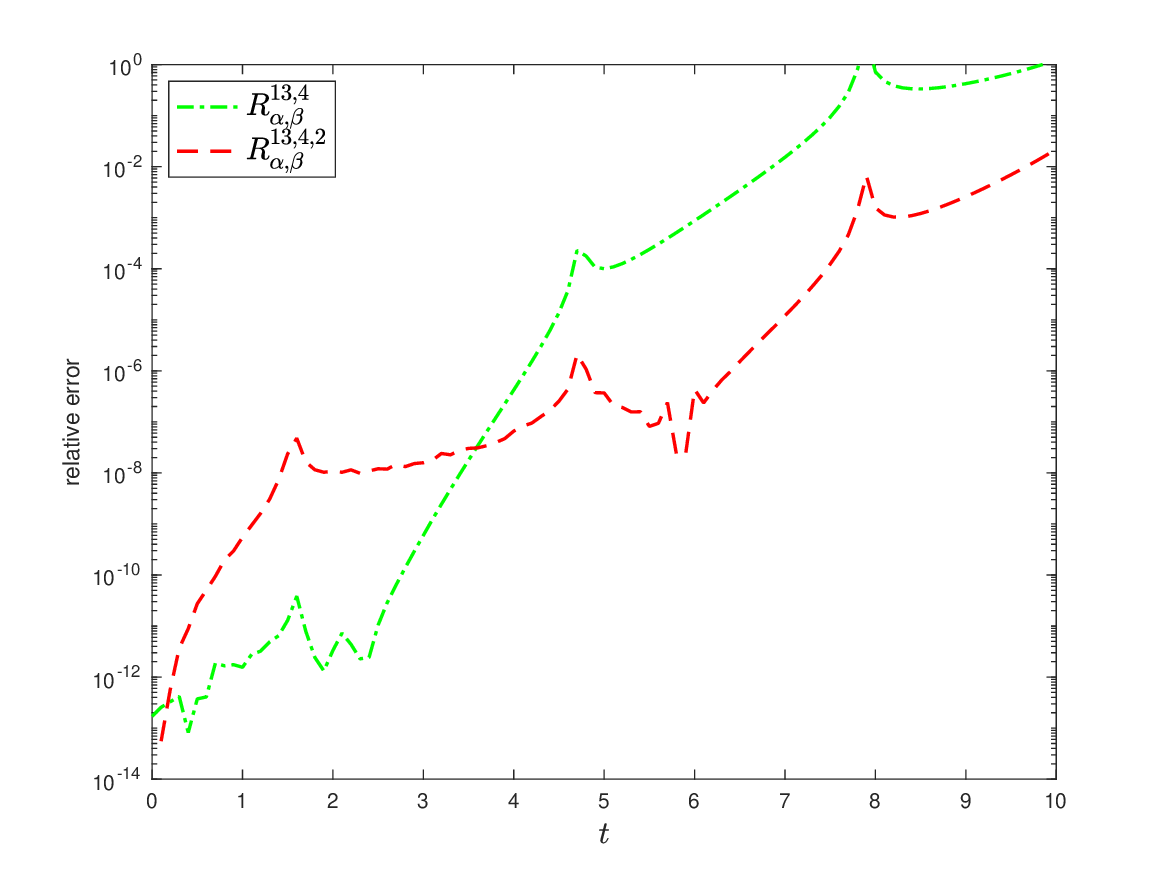}
\caption{Application \ref{wave:ex}: 
approximation of $u(\pi/2,t)$ (left) and relative error 
(right), $\alpha=1.9$.}
\label{fig:wave_a19}
\end{figure}

In Figure \ref{fig:wave_a12}, the solution profile is presented for the case $\alpha = 1.2$. 
As can be seen, both the global Pad\'e approximant $R^{13,4}_{1.2,1}$ and its generalization $R^{13,4,2}_{1.2,1}$ are comparable.
This is expected since the solution has only one root and not much oscillatory behavior is experienced.
In Figure \ref{fig:wave_a19}, the case $\alpha = 1.9$ where the MLF is more oscillatory is presented. 
In this case the rational approximant $R^{13,4}_{1.9,1}$ fails to capture the oscillations for large $t$ and thus yields undesirable approximations. 
However, the modified $R^{13,4,r}_{1.9,1}$ with an
appropriate choice $r = 2$ rectifies this issue and gives more accurate approximate values. 

To illustrate the efficiency of these approximants, their runtime is compared with that of 
ml\_martix Matlab function. 
Table \ref{tab:wave} includes the CPU time for computing the solution  
$U(t) = E_\alpha(-At^\alpha) U_0$ on a time mesh $\{t_n = 0.1n, n = 0, 1, \dots 100\}$ as 
well as the maximum relative error between the approximate and reference values. 
The results show that the approximant $R^{13,4,2}_{1.9,1}$ performs about seven times faster than the 
ml\_matrix. 
Moreover, while the approximant $R^{13,4,2}_{1.9,1}$ on an extended interval is able 
to achieve higher accuracy than $R^{13,4}_{1.9,1}$, their runtime is comparable.


\begin{table}[]
\centering 
\caption {Application \ref{wave:ex}: maximum relative error (RE) and runtime in computing 
the solution of \eqref{eq:wave system} over the time mesh $[0, 0.1, 0.2, \dots, 10]$, 
$\alpha = 1.9 $. }
\label{tab:wave}
\begin{tabular}{lcc}
	& RE    & Runtime \\ \hline
	$R^{13,4}_{1.9, 1}$   & 3.38E+00& 8.00E-02   \\
	$R^{13,4,2}_{1.9, 1}$ & 2.20E-02& 9.95E-02  \\
	ml\_matrix            & -       & 6.796E-01 
\end{tabular}
\end{table}


\section{Concluding Remarks} 

In this paper, a characterization of the oscillatory and monotone behavior of two-parameter Mittag-Leffler function is given. 
Furthermore, generalized rational approximants over extended intervals are developed.
These approximants are based on decomposing Mittag-Leffler function into a root-free one and a polynomial. 
The approximants have good tracking capabilities of the roots and oscillations over
extended intervals. 
Such approximants play an effective role in the implementation of numerical methods for fractional oscillation equations.

\bibliographystyle{IEEEtran}
\bibliography{my_refrences}

\begin{thebibliography}{10}
\providecommand{\url}[1]{#1}
\csname url@samestyle\endcsname
\providecommand{\newblock}{\relax}
\providecommand{\bibinfo}[2]{#2}
\providecommand{\BIBentrySTDinterwordspacing}{\spaceskip=0pt\relax}
\providecommand{\BIBentryALTinterwordstretchfactor}{4}
\providecommand{\BIBentryALTinterwordspacing}{\spaceskip=\fontdimen2\font plus
\BIBentryALTinterwordstretchfactor\fontdimen3\font minus
  \fontdimen4\font\relax}
\providecommand{\BIBforeignlanguage}[2]{{%
\expandafter\ifx\csname l@#1\endcsname\relax
\typeout{** WARNING: IEEEtran.bst: No hyphenation pattern has been}%
\typeout{** loaded for the language `#1'. Using the pattern for}%
\typeout{** the default language instead.}%
\else
\language=\csname l@#1\endcsname
\fi
#2}}
\providecommand{\BIBdecl}{\relax}
\BIBdecl

\bibitem{Kilbas2006}
A.~Kilbas, H.~Srivastava, and J.~Trujullo, \emph{Theory and applications of
  fractional differential equations}.\hskip 1em plus 0.5em minus 0.4em\relax
  Elsevier, 2006.

\bibitem{Gorenflo2014}
R.~Gorenflo, A.~A. Kilbas, F.~Mainardi, and S.~V. Rogosin,
  \emph{Mittag-{L}effler functions. related topics and applications}.\hskip 1em
  plus 0.5em minus 0.4em\relax Springer, 2014.

\bibitem{gorenflo1996fractional}
R.~Gorenflo and F.~Mainardi, \emph{Fractional oscillations and {Mittag-Leffler}
  functions}.\hskip 1em plus 0.5em minus 0.4em\relax Citeseer, 1996.

\bibitem{mainardi1996fractional}
F.~Mainardi, ``Fractional relaxation-oscillation and fractional diffusion-wave
  phenomena,'' \emph{Chaos, Solitons \& Fractals}, vol.~7, no.~9, pp.
  1461--1477, 1996.

\bibitem{achar2001dynamics}
B.~N. Achar, J.~Hanneken, T.~Enck, and T.~Clarke, ``Dynamics of the fractional
  oscillator,'' \emph{Physica A: Statistical Mechanics and its Applications},
  vol. 297, no. 3-4, pp. 361--367, 2001.

\bibitem{stanislavsky2004fractional}
A.~Stanislavsky, ``Fractional oscillator,'' \emph{Physical review E}, vol.~70,
  no.~5, p. 051103, 2004.

\bibitem{stanislavsky2005twist}
A.~A. Stanislavsky, ``{Twist} of fractional oscillations,'' \emph{Physica A:
  Statistical Mechanics and its Applications}, vol. 354, pp. 101--110, 2005.

\bibitem{Gorenflo2002}
R.~Gorenflo, J.~Loutchko, and Y.~Luchko, ``Computation of the
  {M}ittag-{L}effler function {$E_{\alpha,\beta}(z)$} and its derivative,''
  \emph{Fractional Calculus \& Applied Analysis}, vol.~5, no.~4, pp. 491--518,
  2002.

\bibitem{Garrappa2015}
R.~Garrappa, ``Numerical evalution of two and three parameter
  {M}ittag-{L}effler functions,'' \emph{SIAM Journal on Numerical Analysis},
  vol.~53, no.~3, pp. 1350--1369, 2015.

\bibitem{Sarumi2020}
I.~O. Sarumi, K.~M. Furati, and A.~Q.~M. Khaliq, ``Highly accurate global
  {P}adé approximations of generalized {M}ittag–{L}effler function and its
  inverse,'' \emph{Journal of Scientific Computing}, vol.~82, no.~46, 2020.

\bibitem{Sarumi2021}
I.~O. Sarumi, K.~M. Furati, A.~Q.~M. Khaliq, and K.~Mustapha, ``Generalized
  exponential time differencing schemes for stiff fractional systems with
  nonsmooth source term,'' \emph{Journal of Scientific Computing}, vol.~86,
  no.~23, 2021.

\bibitem{Iyiola2018a}
O.~S. Iyiola, E.~O. Asante-Asamani, and B.~A. Wade, ``A real distinct poles
  rational approximation of generalized {M}ittag-{L}effler functions and their
  inverses: applications to fractional calculus,'' \emph{Journal of
  Computational and Applied Mathematics}, vol. 330, pp. 307--317, 2018.

\bibitem{Starovoitov2007}
A.~P. Starovoitov and N.~A. Starovoitova, ``Pad{\'e} approximants of the
  {M}ittag-{L}effler functions,'' \emph{Sbornik Mathematics}, vol. 198, no.~7,
  pp. 1011--1023, 2007.

\bibitem{Borhanifar2015}
A.~Borhanifar and S.~Valizadeh, ``Mittag-{L}effler-{P}ad{\'e} approximations
  for the numerical solution of space and time fractional diffusion
  equations,'' \emph{International Journal of Applied Mathematics Research},
  vol.~4, no.~4, p. 466, 2015.

\bibitem{Winitzki2003}
S.~Winitzki, ``Uniform approximations for transcendental functions,'' in
  \emph{Computational Science and Its Applications --- ICCSA 2003}, V.~Kumar,
  M.~L. Gavrilova, C.~J.~K. Tan, and P.~L'Ecuyer, Eds.\hskip 1em plus 0.5em
  minus 0.4em\relax Berlin, Heidelberg: Springer, 2003, pp. 780--789.

\bibitem{Atkinson2011}
C.~Atkinson and A.~Osseiran, ``Rational solutions for the time-fractional
  diffusion equation,'' \emph{{SIAM} Journal on Applied Mathematics}, vol.~71,
  no.~1, pp. 92--106, 2011.

\bibitem{Zeng2015}
C.~Zeng and Y.~Q. Chen, ``Global {P}ad{\'e} approximations of the generalized
  {M}ittag-{L}effler function and its inverse,'' \emph{Fractional Calculus and
  Applied Analysis}, vol.~18, no.~6, pp. 1492--1506, 2015.

\bibitem{Ingo2017}
C.~Ingo, T.~R. Barrick, A.~G. Webb, and I.~Ronen, ``{A}ccurate {P}ad{\'e}
  global approximations for the {M}ittag-{L}effler function, its inverse, and
  its partial derivatives to efficiently compute convergent power series,''
  \emph{International Journal of Applied and Computational Mathematics},
  vol.~3, no.~2, pp. 347--362, 2017.

\bibitem{hanneken2007enumeration}
J.~W. Hanneken, D.~M. Vaught, and B.~Achar, ``Enumeration of the real zeros of
  the {Mittag}-{Leffler} function {$E_\alpha (z)$}, $1 <\alpha < 2$,'' in
  \emph{Advances in Fractional Calculus}.\hskip 1em plus 0.5em minus
  0.4em\relax Springer, 2007, pp. 15--26.

\bibitem{hanneken2013alpha}
J.~W. Hanneken, B.~Achar, and D.~M. Vaught, ``An alpha-beta phase diagram
  representation of the zeros and properties of the {Mittag-Leffler}
  function,'' \emph{Advances in Mathematical Physics}, vol. 2013, 2013.

\bibitem{duan2014zeros}
J.-S. Duan, Z.~Wang, and S.-Z. Fu, ``The zeros of the solutions of the
  fractional oscillation equation,'' \emph{Fractional Calculus and Applied
  Analysis}, vol.~17, no.~1, pp. 10--22, 2014.

\bibitem{Garrappa2018}
R.~Garrappa and M.~Popolizio, ``Computing the matrix {M}ittag-{L}effler
  function with applications to fractional calculus,'' \emph{Journal of
  Scientific Computing}, vol.~77, no.~1, pp. 129--153, 2018.

\bibitem{honain2023rational}
A.~H. Honain and K.~M. Furati, ``Rational approximation for oscillatory
  mittag-leffler function,'' in \emph{2023 International Conference on
  Fractional Differentiation and Its Applications (ICFDA)}.\hskip 1em plus
  0.5em minus 0.4em\relax IEEE, 2023, pp. 1--5.

\bibitem{haubold2011mittag}
H.~J. Haubold, A.~M. Mathai, and R.~K. Saxena, ``{Mittag-Leffler} functions and
  their applications,'' \emph{Journal of applied mathematics}, vol. 2011, 2011.

\bibitem{Sadeghi2018}
A.~Sadeghi and J.~R. Cardoso, ``Some notes on properties of the matrix
  {M}ittag-{L}effler function,'' \emph{Applied Mathematics and Computation},
  vol. 338, pp. 733--738, 2018.

\bibitem{higham2008functions}
N.~J. Higham, \emph{Functions of matrices: theory and computation}.\hskip 1em
  plus 0.5em minus 0.4em\relax SIAM, 2008.

\bibitem{aguilar2014fractional}
J.~G. Aguilar, J.~R. Hern{\'a}ndez, R.~E. Jim{\'e}nez, C.~Astorga-Zaragoza,
  V.~O. Peregrino, and T.~C. Fraga, ``Fractional electromagnetic waves in
  plasma,'' \emph{Proc. Romanian Acad. A}, vol.~17, no.~1, pp. 31--38, 2014.

\bibitem{murillo2011explicit}
J.~Q. Murillo and S.~B. Yuste, ``{An} explicit difference method for solving
  fractional diffusion and diffusion-wave equations in the {Caputo }form,''
  \emph{Journal of Computational and Nonlinear Dynamics}, vol.~6, no.~2, 2011.

\end{thebibliography}

\end{document}